\documentclass[preprint,12pt]{elsarticle}
\usepackage{mathtools}
\usepackage{lipsum}
\usepackage{amsfonts}
\usepackage{graphicx}
\usepackage{epstopdf}
\usepackage{algorithm}
\usepackage{algorithmic}
\ifpdf
  \DeclareGraphicsExtensions{.eps,.pdf,.png,.jpg}
\else
  \DeclareGraphicsExtensions{.eps}
\fi

\usepackage{adjustbox}
\usepackage{xcolor}
\usepackage{pgfplots}
\usepackage{booktabs}
\usepackage{tikz}
\usetikzlibrary{positioning}
\usetikzlibrary{3d}
\usetikzlibrary{matrix}
\usetikzlibrary{decorations.text}
\usetikzlibrary{spy}

\usepackage{amssymb,amsmath,amsthm}

\makeatletter
\tikzoption{canvas is plane}[]{\@setOxy#1}
\def\@setOxy O(#1,#2,#3)x(#4,#5,#6)y(#7,#8,#9)%
{\def\tikz@plane@origin{\pgfpointxyz{#1}{#2}{#3}}%
	\def\tikz@plane@x{\pgfpointxyz{#4+#1}{#5+#2}{#6+#3}}%
	\def\tikz@plane@y{\pgfpointxyz{#7+#1}{#8+#2}{#9+#3}}%
	\tikz@canvas@is@plane
}
\makeatother
\tikzoption{canvas is xy plane at z}[]{%
	\def\tikz@plane@origin{\pgfpointxyz{0}{0}{#1}}%
	\def\tikz@plane@x{\pgfpointxyz{1}{0}{#1}}%
	\def\tikz@plane@y{\pgfpointxyz{0}{1}{#1}}%
	\tikz@canvas@is@plane
}
\makeatother

\usepackage{subcaption}
\usepackage{stackengine,graphicx,xcolor}
\usepackage{cases}
 
\usepackage{enumitem}

\usepackage{array}
\makeatletter
\newcommand{\thickhline}{%
    \noalign {\ifnum 0=`}\fi \hrule height 1pt
    \futurelet \reserved@a \@xhline
}
\newcolumntype{"}{@{\hskip\tabcolsep\vrule width 1pt\hskip\tabcolsep}}
\makeatother
\newcommand{\TextBB}[1] {\textcolor{blue}{\textbf{#1}} }
\newcommand{\TextGB}[1] {\textcolor{green}{\textbf{#1}} }
\newcommand{\TextB}[1] {\textcolor{blue}{#1} }
\newcommand{\TextG}[1] {\textcolor{green}{#1} }

\theoremstyle{plain}                    %
\newtheorem{theorem}{Theorem}[section]

\theoremstyle{definition}               %
\newtheorem{definition}{Definition}[section]%
\theoremstyle{remark}

\newcommand{\Abold}{\mathbf{A}}

\newcommand{\Dbold}{\mathbf{D}}

\newcommand{\Ibold}{\mathbf{I}}

\newcommand{\Lbold}{\mathbf{L}}
\newcommand{\Mbold}{\mathbf{M}}

 \newcommand{\ebold}{\mathbf{e}}

 \newcommand{\tbold}{\mathbf{t}}
 \newcommand{\ubold}{\mathbf{u}}
 \newcommand{\vbold}{\mathbf{v}}
 \newcommand{\xbold}{\mathbf{x}}

\newcommand{\zbold}{\mathbf{z}}

\usepackage{lineno,hyperref}

\journal{Applied Mathematics and Computation}

\begin{document}

\begin{frontmatter}

\title{Plug-and-Play gradient-based denoisers applied to CT image enhancement}

\author[1]{Pasquale Cascarano}
\author[2]{Elena Loli Piccolomini \corref{cor}}
\ead{elena.loli@unibo.it}
\author[2]{Elena Morotti}
\author[1]{Andrea Sebastiani}

\cortext[cor]{Corresponding author}

\address[1]{Department of Mathematics, University of Bologna, Piazza di Porta S. Donato 5, 40126 Bologna, Italy}
\address[2]{Department of Computer Science and Engineering, University of Bologna, Via Mura Anteo Zamboni 7, 40126 Bologna, Italy}

\begin{abstract}
Blur and noise corrupting Computed Tomography (CT) images can hide or distort small but important details, negatively affecting the diagnosis. In this paper, we present a novel gradient-based Plug-and-Play algorithm, constructed on the Half-Quadratic Splitting scheme, and we apply it to restore CT images. In particular, we consider different schemes encompassing external and internal denoisers as priors, defined on the image gradient domain.
The internal prior is based on the Total Variation functional. The external denoiser is implemented by a deep Convolutional Neural Network (CNN) trained on the gradient domain (and not on the image one, as in state-of-the-art works). We also prove a general fixed-point convergence theorem under weak assumptions on both internal and external denoisers. 
The experiments confirm the effectiveness of the proposed framework in restoring blurred noisy CT images, both in simulated and real medical settings. The achieved enhancements in the restored images are really remarkable, if compared to the results of many state-of-the-art methods. 
\end{abstract}

\begin{keyword}
Plug-and-Play, Half-Quadratic Splitting, external-internal image priors, CNN denoisers, deep learning
\end{keyword}

\end{frontmatter}

\section{Introduction}\label{sec:intro} 
In the field of computational imaging, Image Restoration (IR) aims at recovering an unknown clean image from its noisy and/or blurred measurement.
In Computed Tomography (CT) the presence of blur and noise reduces diagnostic accuracy, hiding or distorting some small but important objects in the reconstructed image. There are different hardware sources of error which cause blur, such as the finite X-ray focal spot size or the spreading effect in the scintillator in Cone Bean Computed Tomography \cite{luo2018cone}. Moreover, quantum noise creating random variations in the attenuation coefficients of X-rays, represents the main contribution to the total noise in CT images. Many statistical analysis have shown that the image noise generated by CT scanner can be regarded as normally distributed \cite{gravel2004method,lu2001noise}. 
Since it is very difficult to avoid these effects by hardware techniques, the software approach is fundamental and several algorithms have been proposed to reduce the blurring and noise artifacts in the CT images.
Examples of restoration algorithms for CT images acquired with different geometries can be found in \cite{al2015deblurring,al2012reducing,jiang2003blind,yim2020deep} and references therein.

Mathematically, by lexicographically reordering the images as vectors, a generic IR task can be written as the following inverse problem:

\begin{equation}\label{eq:inverse_problem}
\text{find} \ \ubold  \quad   \text{such that} \quad  \vbold = \Abold \ubold + \ebold,
\end{equation}
where $\vbold \in \mathbb{R}^{n}$ is the given image, $\ubold \in \mathbb{R}^{n}$ is the unknown desired image and $\Abold \in \mathbb{R}^{n \times n}$ is the forward linear operator defining the IR specific task. The observed image $\vbold$ is usually affected by noise  $\ebold \in \mathbb{R}^{n}$, which  we assume in this work as Additive White Gaussian Noise (AWGN). 

In general, IR problems as \eqref{eq:inverse_problem} are well-known to be ill-posed, meaning that the properties of existence, uniqueness and stability of the desired solution $\ubold$ are not all guaranteed \cite{bertero2020introduction}. Hence, model-based reconstruction methods attempt to find a good estimate $\ubold^{*} \in \mathbb{R}^{n}$ as the solution of a minimization problem whose objective  function is the sum of two  terms $f$ and $g$, namely:  

\begin{equation}\label{eq:problem_f+g}
\ubold^{*} \in \underset{\ubold \in \mathbb{R}^{n}}{\arg\min} \left\{f(\ubold) + g(\ubold)\right\}.
\end{equation}
The functions $f$ and $g$ are usually referred to as \textit{data fidelity} and \textit{regularization} terms, respectively. The former is a task-related term which models the noise affecting the starting measurement $\vbold$, whereas the latter induces prior information on the estimate $\ubold^{*}$ by reflecting, for example,  sparsity patterns, smoothness or geometric assumptions. 
Often, $f$ is set as an L$_{p}$-norm based function measuring the residual between $\Abold \ubold$ and $\vbold$, with $p$  strictly related to noise statistics. 
It is well-known that a squared L$_{2}$-norm fidelity fits with the previous assumption of AWGN affecting the measurement $\vbold$. 

The choice of a regularizer is a challenging task. 
In a model-based approach, a widely used strategy is to define $g$ as a handcrafted term based on desired properties of the reconstructed image in a specific domain, such as the gradient or the wavelet domain, which are effective in medical imaging.  In particular, the Total Variation (TV) \cite{rudin1992nonlinear} is largely employed in this field for its effectiveness in removing noise and preserving curved contours of the objects \cite{landi2012efficient,li2011medical}.\\
These handcrafted regularizers are not completely satisfactory.
However, we believe that a new frontier in the image processing field is represented by the Plug-and-Play (PnP) framework, firstly proposed in \cite{venkatakrishnan2013plug}, where the authors strikingly showed that a closed-form regularizer is not always the best possible choice to properly induce prior information on the desired solution.
Technically, the PnP approach derives from the iterative scheme of proximal algorithms, applied to solve regularized optimization problems as \eqref{eq:problem_f+g}, whose resulting modular structure allows to deal with the data fidelity $f$ and the regularization term $g$, separately. 
Since the sub-step involving $g$ reads as a denoising problem, it can be replaced by any off-the-shelf denoiser, so that the computed solution inherits prior information not necessarily deriving from a closed-form regularization term. 
So far, a large number of papers on PnP have been published focusing on different aspects of the scheme, such as  the proximal algorithm  or the denoiser employed.
Different proximal algorithms are considered such as the Alternating Direction Method of Multipliers (ADMM), the Half-Quadratic Splitting (HQS) or the Fast Iterative Shrinkage-Thresholding Algorithm (FISTA) \cite{venkatakrishnan2013plug,sreehari2016plug, zhang2017learning,kamilov2017plug}. \\
Moreover, several denoisers have already been successfully used and they are usually labelled as internal or external denoisers \cite{mosseri2013}. 

Internal denoisers, inducing \textit{internal priors}, are tailored to define features onto the observed data but they struggle to deal with several different image features simultaneously. Examples are the proximal maps of handcrafted regularizers, the BM3D \cite{dabov2007image} and the Non-Local Mean (NLM) filter \cite{buades2005non}. 
External denoisers, inducing \textit{external priors}, are related to an outer set of clean images, so they can fail when dealing with unseen noise variance and image patterns.
Early studies made use of Gaussian Mixture Models (GMMs) \cite{zoran2011learning} and trained nonlinear reaction diffusion based denoisers \cite{chen2016trainable} as external denoisers.
Since nowadays deep learning based priors lead to outstanding performances for denoising images \cite{xie2012image,burger2012image},  PnP frameworks are equipped with  pre-trained Convolutional Neural Network (CNN) denoisers  as in \cite{zhang2017learning,zhang2020plug,meinhardt2017learning}. 

The aforementioned approaches exploit either external or inetrnal denoisers; very recently some generalizations to handle multiple internal and/or external denoisers have been proposed in \cite{rond2016poisson,he2021support}.

\subsection*{Motivation and contributions of the paper}
In all the previously aforementioned works on PnP, the deep learning based denoisers act on the image domain. Since it is well-known that a prior defined on the gradient domain may enhance the image reconstructions both in terms of shape recovering and noise removal, we propose here a PnP method specifying a gradient-based external prior through CNN networks trained to restore the corrupted image gradients. We refer to this method as GCNN, in the following. 
Moreover, motivated by the apparent complementarity of external and internal denoisers, we propose a new hybrid PnP framework, relying on the Half-Quadratic Splitting algorithm, which combines the Total Variation and a CNN-based denoiser, acting either on the image domain or on the image gradient domain. The former and the latter are referred as ICNN-TV and GCNN-TV, respectively. \\
Finally, a fixed point convergence is proved for all the implemented methods upon weak assumptions on the considered denoisers.

We test GCNN, ICNN-TV and GCNN-TV to restore blurred and noisy synthetic and real CT images.
The performances of our proposals are validated through several comparisons with other state-of-the-art PnP methods exploiting different denoisers. 
The numerical results provide very high quality reconstructions and confirm the robustness of the proposed gradient-based frameworks both in restoring low contrasted objects with different shapes and in removing noise.

\subsection*{Organization of the paper}
In this paper, we present in Section \ref{sec:method} the proposed PnP methods together with some implementation choices for the denoisers implemented. In Section \ref{sec:experiments} we report and analyse the numerical results. Finally, in Section \ref{sec:conclusions}, we conclude the paper with a brief discussion. In Appendix \ref{sec:theorem} we report a fixed-point convergence theorem for the proposed schemes and its proof.

\section{Proposed Plug-and-Play methods} \label{sec:method}

We describe here the proposed algorithms for the solution of problem \eqref{eq:problem_f+g} with particular choices of $f$ and $g$.

Due to the previous assumption of AWGN affecting the measurement $\vbold$, we fix the fidelity term as $f(\ubold):=\frac{1}{2}\lVert \Abold \ubold - \vbold \rVert^{2}_{2}$.
Concerning the regularizer $g$
we propose to set it as the sum of two terms $g_1$ and $g_2$, weighted by the nonnegative parameters  $\lambda$ and $\eta$, respectively. 
Hence, the minimization problem \eqref{eq:problem_f+g} reads:

\begin{equation} \label{eq:l2+g}
\ubold^{*} \in\underset{\ubold \in \mathbb{R}^{n}}{\arg\min}\left\{\frac{1}{2}\lVert \Abold\ubold-\vbold\rVert_2^2 + \lambda g_{1}(\Lbold_1\ubold)+\eta g_{2}(\Lbold_2\ubold) \right\}.
\end{equation}
Moreover, we assume $g_{1}$ and $g_{2}$ are positive and convex real-valued maps:

\begin{equation}
g_{1}:\mathbb{R}^{l_1}\rightarrow\mathbb{R}^+, \quad g_{2}:\mathbb{R}^{l_2}\rightarrow\mathbb{R}^+,
\end{equation}
with $l_1$ and $l_2$ positive integers, $\Lbold_1 \in \mathbb{R}^{l_1 \times n}$ and $\Lbold_2 \in \mathbb{R}^{l_2 \times n}$. 
 
We now consider the HQS iterative method described in \cite{geman1995nonlinear,wang2008new} as numerical solver to compute $\ubold^{*}$. By introducing the auxiliary variables $\tbold \in \mathbb{R}^{l_{1}}$ and $\zbold \in \mathbb{R}^{l_{2}}$ subject to $\tbold := \Lbold_1 \ubold$ and $\zbold := \Lbold_2 \ubold$, the following penalized half-quadratic function is taken into account: 

\begin{equation}
\begin{split}
\mathcal{L}(\ubold, \tbold, \zbold; \rho^{\tbold}, \rho^{\zbold}):= & \ 
\frac{1}{2}\lVert \Abold\ubold-\vbold\rVert_2^2 + \lambda g_{1}(\tbold)+\eta g_{2}(\zbold) +  \\ 
 & +  \frac{\rho^{\tbold}}{2}\lVert \Lbold_1 \ubold-\tbold\rVert_2^2 + \frac{\rho^{\zbold}}{2}\lVert \Lbold_2 \ubold-\zbold\rVert_2^2. 
\end{split}  
\label{eq:cost_function}
\end{equation}

At each iteration $k$, the HQS algorithm performs this alternated minimization scheme with respect to $\tbold$, $\zbold$ and the primal variable $\ubold$:  

\begin{numcases}{}
\scalebox{0.93}{$\tbold_{k+1} \in \underset{\tbold\in\mathbb{R}^{l_1}}{\arg\min}~ \lambda g_{1}(\tbold)+\dfrac{\rho^{\tbold}_k}{2}\lVert \Lbold_1 \ubold_k - \tbold\rVert_2^2$} & \label{eq:sub_t2} \\
\scalebox{0.93}{$\zbold_{k+1} \in \underset{\zbold\in\mathbb{R}^{l_2}}{\arg\min}~ \eta g_{2}(\zbold)+\dfrac{\rho^{\zbold}_k}{2}\lVert \Lbold_2 \ubold_k - \zbold\rVert_2^2$} \label{eq:sub_z2}& \\
\begin{aligned}
\scalebox{0.93}{$\ubold_{k+1}=\underset{\ubold \in \mathbb{R}^{n}}{\arg\min} ~ \dfrac{1}{2}\lVert \Abold\ubold-\vbold\rVert_2^2 + \dfrac{\rho^{\tbold}_k}{2}\lVert \Lbold_1 \ubold - \tbold_{k+1}\rVert_2^2 + \dfrac{\rho^{\zbold}_k}{2}\lVert \Lbold_2 \ubold - \zbold_{k+1}\rVert_2^2$}, 
\end{aligned}
\label{eq:sub_u2} & 
\end{numcases}

where $(\rho^{\tbold}_k)_{k=1}^{\infty}$ and $(\rho^{\zbold}_k)_{k=1}^{\infty}$ are two non-decreasing sequences of positive penalty parameters.
The key feature of HQS is that the prior related sub-steps \eqref{eq:sub_t2} and \eqref{eq:sub_z2} are specified through the proximal maps of $g_{1}$ and $g_{2}$, respectively, which are mathematically equivalent to regularized denoising problems. The PnP framework exploits both this equivalence and the modular structure of the algorithm  by replacing such proximal maps with any off-the-shelf denoiser.

Now we present the proposed hybrid PnP scheme. We introduce a pre-trained learning-based denoiser $\mathcal{D}_{\sigma}^{\text{ext}}$ and an image-specific denoiser $\mathcal{D}_{\gamma}^{\text{int}}$ depending on the positive parameters $\sigma$ and $\gamma$ related to the noise-level in the images to be denoised, so that the greater are $\sigma$ and $\gamma$, the more powerful is the denoising realized. In particular, in our scheme we choose two sequences $(\sigma_{k})_{k=1}^{+ \infty}$ and $(\gamma_{k})_{k=1}^{+ \infty}$ such that, at step $k$, $\mathcal{D}_{\sigma_k}^{\text{ext}}$ and $\mathcal{D}_{\gamma_k}^{\text{int}}$ replace the sub-steps \eqref{eq:sub_t2} and \eqref{eq:sub_z2}, respectively. A standard assumption in PnP is that $\sigma_{k}$ and $\gamma_{k}$ are both related with the penalty parameters $\rho_{k}^{\tbold}$ and $\rho_{k}^{\zbold}$ through these formulas:
\begin{equation}
    \sigma_k:= \sqrt{\frac{\alpha}{{\rho_{k}^{\tbold}}}}, \ \ \gamma_k:= \sqrt{\frac{\beta}{{\rho_{k}^{\zbold}}}},
    \label{eq:sigma_gamma}
\end{equation}  where $\alpha$ and $\beta$ are chosen positive scaling factors. A sketch  of the resulting hybrid PnP framework is reported in Algorithm \ref{alg:hybrid_pnp_hqs}.

\begin{algorithm}
\caption{}\label{alg:hybrid_pnp_hqs} %
\begin{algorithmic}
	\STATE \textbf{Input:} $\alpha$, $\beta$ and  $(\rho^{\tbold}_k)_{k=1}^{\infty}$,\   $(\rho^{\zbold}_k)_{k=1}^{\infty}$, $\Abold$, $\Lbold_1$, $\Lbold_2$, $\vbold$, $\ubold_{1}$, \textit{K}.
	\FOR{k = 1 \dots \textit{K}} 
	\STATE $\tbold_{k+1}=\mathcal{D}^{\text{ext}}_{\sigma_k}(\Lbold_1 \ubold_k) $
	\STATE $\zbold_{k+1}=\mathcal{D}^{\text{int}}_{\gamma_k}(\Lbold_2 \ubold_k)$
	\STATE $\ubold_{k+1}= \underset{\ubold \in \mathbb{R}^{n}}{\arg\min}\dfrac{1}{2}\lVert \Abold\ubold-\vbold\rVert_2^2+\dfrac{\rho^{\tbold}_k}{2}\lVert \Lbold_1 \ubold - \tbold_{k+1}\rVert_2^2 + \dfrac{\rho^{\zbold}_k}{2}\lVert \Lbold_2 \ubold - \zbold_{k+1}\rVert_2^2$
	\ENDFOR
\end{algorithmic}
\end{algorithm} 

Under some quite general assumptions on the denoisers and on the sequences $\rho_{k}^{\tbold}$ and $\rho_{k}^{\zbold}$, the iterates defined in Algorithm \ref{alg:hybrid_pnp_hqs} converge to a fixed-point $(u^*,t^*,z^*)$. See Appendix \ref{sec:theorem}. for an in-depth discussion on the hypothesis and the fixed-point convergence theorem.

We describe now the choices of $\mathcal{D}_{\sigma}^{\text{ext}}$ and $\mathcal{D}_{\sigma}^{\text{int}}$ in our particular implementation.
We fix as internal denoiser a scheme based on the Total Variation (TV) \cite{rudin1992nonlinear}. The properties of edge preserving and noise suppressing of the TV in many image processing applications are well-established.  The TV function is defined as:  
\begin{equation}   \label{eq:TV_reg}
\text{TV}(\ubold) :=  \sum_{i=1}^{n} \lVert (\Dbold \ubold)_{i} \rVert_{2}  = \sum_{i=1}^{n} \left( \sqrt{(\Dbold_h \ubold)_i^2 + (\Dbold_v \ubold)_i^2}\right),
\end{equation}
where $(\Dbold \ubold)_i:=((\Dbold_{h} \ubold)_i,(\Dbold_{v} \ubold)_i)  \in \mathbb{R}^{2}$, for $i=1 \dots n$  denotes the discrete image gradient computed at pixel $i$ along the horizontal and vertical axes, separately. 
Hence,  the function $g_{2}$ in \eqref{eq:l2+g} is set as:
\begin{align} \label{eq:g2}
    g_{2}: \quad  & \mathbb{R}^{2 \times n} \to \mathbb{R} \nonumber \\
                    &   \xbold      \to    \sum_{i=1}^{n} \lVert \xbold_{i} \rVert_{2} \quad \text{with} \quad \xbold_i \in \mathbb{R}^{2},
\end{align}
and the $\Lbold_{2}$ operator is the discrete gradient $\Dbold = (\Dbold_{h};\Dbold_{v})$ such that $\Dbold_{h},\Dbold_{v} \in \mathbb{R}^{n \times n}$ are the finite differences discretization of first order derivative operators along the horizontal and vertical axes, respectively.
We remark that, in Algorithm \ref{alg:hybrid_pnp_hqs}, $\mathcal{D}^{\text{int}}_{\gamma_k}$ is the proximal map of $g_{2}$ with parameter $\gamma_{k}^{2}=\frac{\eta}{\rho_{k}^{\zbold}}$. 

Concerning the choice of the external denoiser, due to the state-of-the-art performances in denoising tasks reached by deep learning strategies \cite{xie2012image,burger2012image}, we embed a Deep CNN denoiser $\mathcal{D}_{\sigma}^{\text{CNN}}$ in our Algorithm \ref{alg:hybrid_pnp_hqs} as $\mathcal{D}_{\sigma}^{\text{ext}}$. Previous studies have already successfully inspected a CNN-based PnP \cite{zhang2017learning,meinhardt2017learning} whose CNN denoisers act directly only on the image-domain. Conversely, our denoiser acts on the image through an operator $\Lbold_1$, which we set either equal to the identity matrix $\mathbf{I}$ or the discrete derivative $\Dbold$. We will explain in \ref{subsec:impl} how we have implemented the action of the CNN with respect to the choice of the operator  $\Lbold_1$.

The method obtained with the described choices of CNN as external denoiser and TV functional as internal denoiser is reported in Algorithm \ref{alg:pnp_hqsTv}. In the following, we will denote it as  ICNN-TV or GCNN-TV when $\Lbold_{1}=\mathbf{I}$ or $\Lbold_{1}=\Dbold$, respectively.

\begin{algorithm}
\caption{}\label{alg:pnp_hqsTv} %
\begin{algorithmic}
	\STATE \textbf{Input:} $\alpha$, $\beta$ and  $(\rho^{\tbold}_k)_{k=1}^{\infty}$,\   $(\rho^{\zbold}_k)_{k=1}^{\infty}$, $\Abold$, $\Lbold_1$, $\vbold$, $\ubold_{1},$ \textit{K}.
	\FOR{k = 1 \dots \textit{K}} 
	\STATE $\tbold_{k+1}=\mathcal{D}^{\text{CNN}}_{\sigma_k}(\Lbold_1 \ubold_k)$
	\STATE $\zbold_{k+1}= \text{prox}_{g_{2}} (\Dbold \ubold_k)$
    \STATE $\ubold_{k+1}=\underset{\ubold \in \mathbb{R}^{n}}{\arg\min}\dfrac{1}{2}\lVert \Abold\ubold-\vbold\rVert_2^2+\dfrac{\rho^{\tbold}_k}{2}\lVert \Lbold_1 \ubold - \tbold_{k+1}\rVert_2^2 + \dfrac{\rho^{\zbold}_k}{2}\lVert \Lbold_2 \ubold - \zbold_{k+1}\rVert_2^2$
	\ENDFOR
\end{algorithmic}
\end{algorithm} 

We also consider the case where only the external CNN denoiser is plugged in (thus excluding the internal prior), to better investigate the effectiveness provided by the use of the sole gradient-based CNN. We report the general scheme in Algorithm \ref{alg:pnp_hqs}. 
Coherently, we label this scheme as ICNN or GCNN when $\Lbold_{1}=\mathbf{I}$ or $\Lbold_{1}=\Dbold$, respectively. We stress that ICNN is equivalent to the approach proposed in \cite{zhang2017learning} and there named as IRCNN.

\begin{algorithm}
\caption{} \label{alg:pnp_hqs} %
\begin{algorithmic}
	\STATE \textbf{Input:} $\alpha$ and  $(\rho^{\tbold}_k)_{k=1}^{\infty}$,\ $\Abold$, $\Lbold_1$, $\vbold$, $\ubold_{1},$ \textit{K}.
	\FOR{k = 1 \dots \textit{K}} 
	\STATE $\tbold_{k+1}=\mathcal{D}^{\text{CNN}}_{\sigma_k}(\Lbold_{1} \ubold_k)$
	\STATE $\ubold_{k+1}=\underset{\ubold}{\arg\min}\left\{\dfrac{1}{2}\lVert \Abold\ubold-\vbold\rVert_2^2+\dfrac{\rho^{\tbold}_k}{2}\lVert \Lbold_{1}\ubold - \tbold_{k+1}\rVert_2^2 \right\}$
	\ENDFOR
 \end{algorithmic}
 \end{algorithm} 

\subsection{Implementation notes \label{subsec:impl}}

We now refer to particular implementation choices when the proposed Algorithm \ref{alg:pnp_hqsTv} and Algorithm \ref{alg:pnp_hqs} are applied to image deblur, as considered in our numerical experiments.
In Algorithm \ref{alg:pnp_hqsTv}, at each iteration $k$, the minimization problem on the primal variable $\ubold$ is solved by applying the first order optimality conditions leading to the following linear system:
\begin{equation}\label{eq:normeq}
    (\Abold^T\Abold+\rho^{\tbold}_k\Lbold_1^T\Lbold_1+\rho^{\zbold}_k \Dbold^T\Dbold)\ubold_{k+1}=\Abold^T\vbold+\rho^{\tbold}_k\Lbold_1^T\tbold_{k+1}+\rho^{\zbold}_k \Dbold^T\zbold_{k+1}.
\end{equation}
The matrix $\Abold$ is the Toeplitz matrix corresponding to the discrete two dimensional convolution with a specific Point Spread Function (PSF). If we assume periodic boundary conditions, all the previous linear operators are diagonalized by the Fourier Transform and the linear system \eqref{eq:normeq} can be efficiently solved by using the Fast Fourier Transform \cite{hansen2008deblurring}. Similarly, $\ubold_{k+1}$ in Algorithm \ref{alg:pnp_hqs} can be computed by setting $\rho^{\zbold}_k=0$ in \eqref{eq:normeq}.  

Concerning  the update of $\zbold_{k}$ in Algorithm \ref{alg:pnp_hqsTv}, we observe that it reduces to the solution of  $n$ bi-dimensional optimization problems which can be computed in a closed form by using the proximal map of the L$_2$-norm.

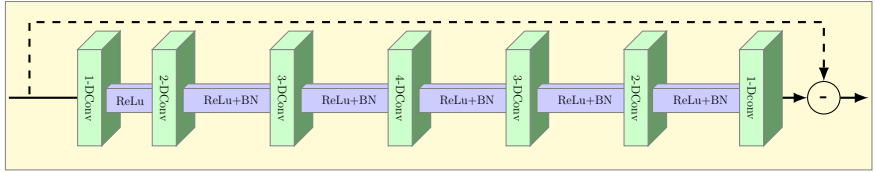
\begin{figure*}[!ht]
						\centering
						\begin{tikzpicture}[scale=0.32, transform shape]
						\tikzset{pics/fake box/.style args={%
								#1 with dimensions #2 and #3 and #4 rot #5 text #6}{
								code={
									\draw[gray,ultra thin,fill=#1!20!gray] (0,0,0) -- ++(0,0,-#4) -- ++(0,-#3,0) -- ++(0,0,#4) -- cycle;
									\draw[gray,ultra thin,fill=#1!20] (0,0,0) -- ++(-#2,0,0) -- ++(0,0,-#4)-- ++(#2,0,0) -- cycle;
									\draw[gray,ultra thin,fill=#1!20] (0,0,0) -- ++(-#2,0,0) -- ++(0,-#3,0) -- ++(#2,0,0) -- cycle;
									\node at (-#2/2,-#3/2) [rotate=#5, anchor=center] {#6};
								}
								
						}}
						\node[] (p00) at (-3,-2){};
						\node[] (p01) at (-4,-2){};
						\draw[gray,ultra thin,fill=yellow!20] (-4,2) rectangle (32,-5);
						\node[circle,draw,anchor=center, minimum size=1.35cm] (p1) at (30,-2,0) {\Huge{\textbf{-}}};
						\draw[thick,-latex] (p01.east) -- (p1);
						\draw pic (layer0) at (0,0,0) {fake box={green} with dimensions 1 and 4 and 2 rot -90 text {\large 1-DConv}};
						\draw pic (relu1) at (1.8,-2,-1) {fake box={blue} with dimensions 2 and 1 and 0.5 rot 0 text {\large ReLu}};
						\draw pic (layer1) at (3.1,0,0) {fake box={green} with dimensions 1 and 4 and 2 rot -90 text {\large 2-DConv}};
						\draw pic (relu) at (7,-2,-1) {fake box={blue} with dimensions 4 and 1 and 0.5 rot 0 text {\large ReLu+BN}};
						\draw pic (layer2) at (8,0,0) {fake box={green} with dimensions 1 and 4 and 2 rot -90 text {\large 3-DConv}};
						\draw pic (relu) at (11.9,-2,-1) {fake box={blue} with dimensions 4 and 1 and 0.5 rot 0 text {\large ReLu+BN}};
						\draw pic (layer3) at (12.9,0,0) {fake box={green} with dimensions 1 and 4 and 2 rot -90 text {\large 4-DConv}};
						\draw pic (relu) at (16.8,-2,-1) {fake box={blue} with dimensions 4 and 1 and 0.5 rot 0 text {\large ReLu+BN}};
						\draw pic (layer4) at (17.8,0,0) {fake box={green} with dimensions 1 and 4 and 2 rot -90 text {\large 3-DConv}};
						\draw pic (relu) at (21.7,-2,-1) {fake box={blue} with dimensions 4 and 1 and 0.5 rot 0 text {\large ReLu+BN}};
						\draw pic (layer5) at (22.7,0,0) {fake box={green} with dimensions 1 and 4 and 2 rot -90 text {\large 2-DConv}};
						\draw pic (relu) at (26.5,-2,-1) {fake box={blue} with dimensions 4 and 1 and 0.5 rot 0 text {\large ReLu+BN}};
						\draw pic (layer6) at (27.5,0,0) {fake box={green} with dimensions 1 and 4 and 2 rot -90 text {\large 1-Dconv}};
						\draw[thick, dashed,-latex] (p00.north) -- ++(0,3,0) -| (p1);
						
						\node[] (p2) at (32,-2)
						{};%
						\draw[thick,-latex] (p1) -- (p2);
						\end{tikzpicture}
						\caption{DnCNN architecture scheme \cite{zhang2017learning}. BN represents the batch normalization and $m$-DConv denotes $m$-dilated convolution.}\label{fig:ircnn_arch}
					\end{figure*}

\begin{figure*}[!ht]
						\begin{subfigure}{0.47\textwidth}
						\centering
						\begin{tikzpicture}[scale=0.32, transform shape]
						\tikzset{pics/fake box/.style args={%
								#1 with dimensions #2 and #3 and #4 rot #5 text #6}{
								code={
									\draw[gray,ultra thin,fill=#1!20!gray] (0,0,0) -- ++(0,0,-#4) -- ++(0,-#3,0) -- ++(0,0,#4) -- cycle;
									\draw[gray,ultra thin,fill=#1!20] (0,0,0) -- ++(-#2,0,0) -- ++(0,0,-#4)-- ++(#2,0,0) -- cycle;
									\draw[gray,ultra thin,fill=#1!20] (0,0,0) -- ++(-#2,0,0) -- ++(0,-#3,0) -- ++(#2,0,0) -- cycle;
									\node at (-#2/2,-#3/2) [rotate=#5, anchor=center] {#6};
								}
								
						}}
						\node[] (p0) at (-6,-2){\includegraphics[width=5.5cm]{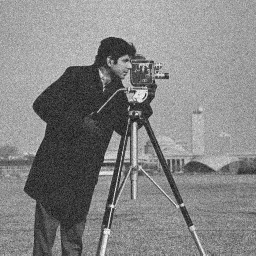}};
						\node[] (p1) at (5,-2){\includegraphics[width=5.5cm]{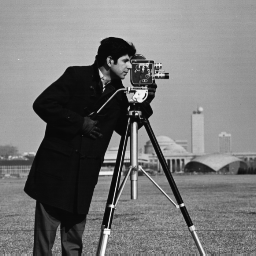}};
						\draw[thick,-latex] (p0.east) -- (p1.west);
						\draw pic (layer0) at (0,0,0) {fake box={yellow} with dimensions 2 and 4.5 and 2 rot -90 text {\Huge DnCNN}};
						\end{tikzpicture}
						\caption{CNN image denoiser on the image.}
						\label{fig:ircnn_applied1}
						\end{subfigure} \begin{subfigure}{0.47\textwidth}
						\centering
						\begin{tikzpicture}[scale=0.32, transform shape]
						\tikzset{pics/fake box/.style args={%
								#1 with dimensions #2 and #3 and #4 rot #5 text #6}{
								code={
									\draw[gray,ultra thin,fill=#1!20!gray] (0,0,0) -- ++(0,0,-#4) -- ++(0,-#3,0) -- ++(0,0,#4) -- cycle;
									\draw[gray,ultra thin,fill=#1!20] (0,0,0) -- ++(-#2,0,0) -- ++(0,0,-#4)-- ++(#2,0,0) -- cycle;
									\draw[gray,ultra thin,fill=#1!20] (0,0,0) -- ++(-#2,0,0) -- ++(0,-#3,0) -- ++(#2,0,0) -- cycle;
									\node at (-#2/2,-#3/2) [rotate=#5, anchor=center] {#6};
								}
								
						}}
						\node[] (p0) at (-6,-2){\includegraphics[width=5.5cm]{architecture/noisy.png}};
						\node[] (p1) at (9.5,-2){\includegraphics[width=5.5cm]{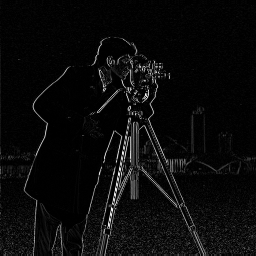}};
						\draw[thick,-latex] (p0.east) -- (p1.west);
						\draw pic (layer0) at (0,0,0) {fake box={yellow} with dimensions 2 and 4.5 and 2 rot -90 text {\Huge DnCNN}};
						\draw pic (layer1) at (4,0,0) {fake box={green} with dimensions 1.5 and 4.5 and 2 rot -90 text  {\Huge FE}};
						\end{tikzpicture}
						\caption{CNN denoiser for gradient restoration.}
						\label{fig:ircnn_applied2}
						\end{subfigure}

						\caption {Trained schemes for denoising.}
						\label{fig:ircnn_applied}
						
					\end{figure*}

To implement the CNN based external denoiser $\mathcal{D}^{\text{CNN}}_{\sigma}$ we adopt the widely used DnCNN architecture proposed in \cite{zhang2017learning} and depicted in Figure \ref{fig:ircnn_arch}. It is constituted by seven dilated convolutional layers \cite{yu2015multi} activated by ReLu functions. 
 We consider for the CNN training  the \textit{Train400} image dataset \cite{chen2016trainable}, also used in \cite{zhang2017learning}. It contains $400$ gray-scale natural images of size $180\times 180$ obtained by cropping larger images in the Berkeley Segmentation dataset \cite{MartinFTM01}.
The training phase of the network depends on the setting of $\Lbold_1$ which  defines the image domain where the noise should be removed. In particular  when we set $\Lbold_1=\Ibold$, the CNN denoiser is trained on noisy images (Figure \ref{fig:ircnn_applied1}) and  we make use of the 25 denoisers downloaded from \url{https://github.com/cszn/IRCNN}, each one trained on a single noise level in the range [2, 50] with step 2. 
Our proposal considers the case $\Lbold_1=\Dbold$. Therefore, the CNN is trained to give as output the noisy-free gradient of the $40\times 40$ patches (Figure \ref{fig:ircnn_applied2}) corrupted with the same 25 noise levels. 
In particular, as depicted in Figure \ref{fig:ircnn_applied2}, we add a general linear Feature Extractor (FE, green layer) at the end of the network, which is taken in our implementation equals to the discrete image gradient operator. %
We use the ADAM optimizer with the Tensorflow default parameters and we set the epochs number to 150. The correspondence between the iteration $k$ of the algorithms and one of the 25 available networks is performed as in \cite{zhang2017learning}.

\section{Results and discussion}\label{sec:experiments}

In this section, we describe the results obtained by testing our GCNN, ICNN-TV and GCNN-TV schemes on the task of image denoising and deblurring. The Python codes of our proposals are available at \url{https://github.com/sedaboni/PnP-TV}. 
We  validate our methods both on a synthetic image, characterized by elements of interest for CT medical purposes, and on two real CT images. %
All the ground truth images have values in the range $[0, 255]$.

Our methods are  compared with the baseline TV regularization implemented in the standard ADMM algorithm, which uses the discrepancy principle \cite{morozov1984methods} for the estimation of the regularization parameter, the IRCNN approach \cite{zhang2017learning} which is referred to as ICNN in the following, the standard PnP with BM3D and NLM chosen as denoisers and a very recent method \cite{he2021support} which combines a truncated L$_{1}$-norm computed on the wavelet operator applied to the signal and BM3D, therefore two internal denoisers.

For a quality assessment of the results, we create artificially blurred and noisy images from a ground truth (GT) image and we compute the Structural Similarity Index Measure (SSIM) and the Peak Signal-to-Noise-Ratio (PSNR) \cite{hore2010image} between the restored image and the ground-truth. 
Moreover, to quantify  noise removal, we compute the standard deviation on uniform regions of interest of the restored images.

For all the proposed algorithms the input parameters $\alpha$ and $\beta$  are heuristically chosen to compute a solution satisfying the discrepancy principle. The algorithm performs at most  30 iterations. The first iterate $\ubold_{1}$ is initialized as a vector of zeros. Concerning the choice of $(\rho^{\tbold}_k)_{k=1}^{\infty}$ and $(\rho^{\zbold}_k)_{k=1}^{\infty}$, we have set $\rho^{\tbold}_k=\rho^{\zbold}_k=k(1 + \epsilon)^{k}$, with $\epsilon >0$, satisfying the conditions required in the fixed-point convergence theorem stated in Appendix \ref{sec:theorem}. All the hyperparameters of the competitors have been fixed in order to provide a solution which satisfies the discrepancy principle. 

\subsection{Results on a synthetic test problem}\label{subsec:ImageDenoisingDeblurring} %
\begin{figure}
	\centering
	\begin{subfigure}[b]{0.49\textwidth}
	\centering
	\begin{tikzpicture}
	\begin{scope}[spy using outlines={rectangle,red,magnification=2,size=1.5cm}]
	\node [name=c]{\frame{\includegraphics[height=4cm]{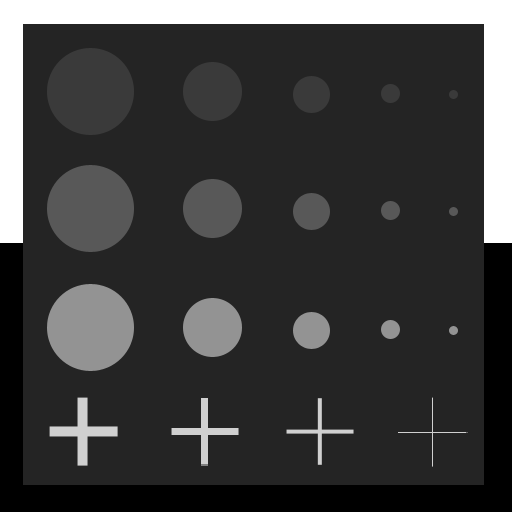}}};
	\draw[green,line width=0.1mm] (-0.15,0.55) rectangle (0.35, 1.05);
	\spy [magnification=1.5] on (-1.28,1.28) in node [name=c1]  at (3.4,1.5);
	\spy[magnification=2] on (1.28,1.28) in node [name=c1]  at (3.4,-0.05);
	\spy [magnification=2] on (1.38,-1.37) in node [name=c1]  at (3.4,-1.6);
	\end{scope}
	\node[] at (2.25,0) {
	\includegraphics[trim=22cm 0.7cm 0cm 0.4cm, clip, height=4cm]{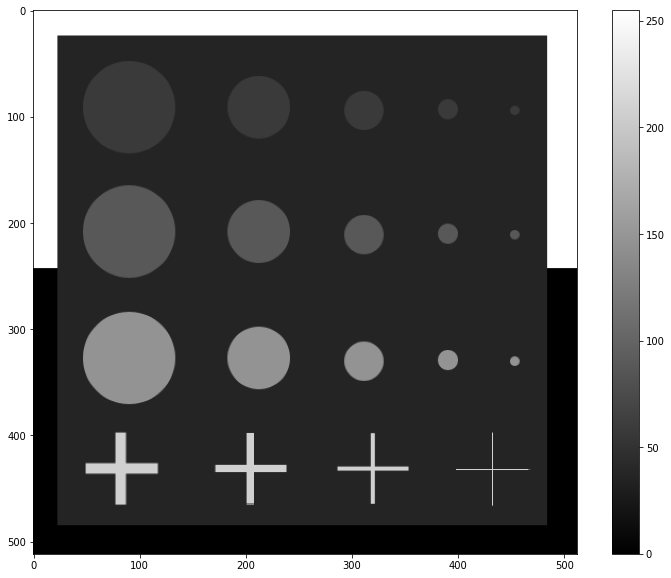}
	};
	\end{tikzpicture}
	\caption{Ground truth}
	\label{fig:cirs-A}
	\end{subfigure}
	\begin{subfigure}[b]{0.49\textwidth}
	\centering
	\begin{tikzpicture}
	\begin{scope}[spy using outlines={rectangle,red,magnification=2,size=1.5cm}]
	\node [name=c]{\frame{\includegraphics[height=4cm]{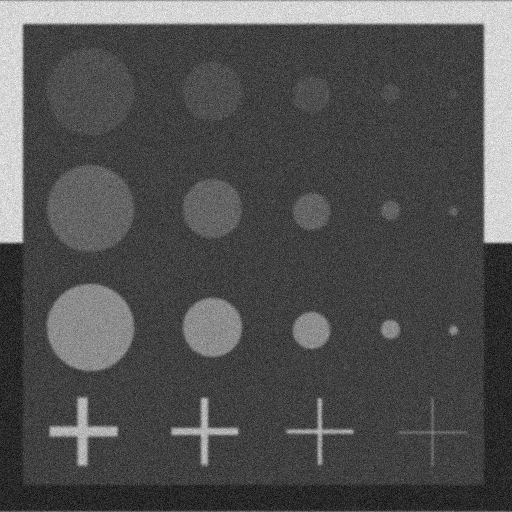}}};
	\spy [magnification=1.5] on (-1.28,1.28) in node [name=c1]  at (3.4,1.5);
	\spy[magnification=2] on (1.28,1.28) in node [name=c1]  at (3.4,-0.05);
	\spy [magnification=2] on (1.38,-1.37) in node [name=c1]  at (3.4,-1.6);
	\end{scope}
	\end{tikzpicture}
	\caption{Corrupted}
	\label{fig:cirs-B}
	\end{subfigure}
    \caption{Ground truth gray-scale test image and a simulated degraded acquisition. In (a) the green square highlights the uniform patch used to evaluate ROI-std. In (a) and (b) three close-ups (red boxes) are depicted alongside the images.}
     \label{fig:cirs}
    \end{figure}

We start our experiments by considering the numerical simulation acting on the gray-scale $512\times 512$ synthetic image reported in Figure \ref{fig:cirs-A}. 
The  image  is designed to test the algorithms performances in the case of low and high contrast objects, with curved and straight borders: the ground truth image contains many circles of different diameter but uniform intensity; each row has homogeneous circles, with enhancing contrast with respect to the uniform background. The fourth row contains crosses of different thickness and high contrast. 
To build our test problems, we blur the ground truth image using a Gaussian $15\times 15$ kernel with zero mean and standard deviation $1.2$, then we introduce AWGN with standard deviation $std$ in $\{10, 15, 20 \}$.
In Figure \ref{fig:cirs-B} we show the corrupted image obtained with $std=15$. In Figure \ref{fig:cirs-A} and \ref{fig:cirs-B},  we also depict three close-ups on the regions bounded by red squares.  

\begin{figure}
	\centering
	\scalebox{0.65}{
    \begin{subfigure}[b]{0.18\textwidth}
	\centering \begin{tikzpicture}
	\begin{scope}[spy using outlines={rectangle,red,magnification=4,size=0.9\textwidth}]
	\node [name=c]{{\includegraphics[height=0.9\textwidth]{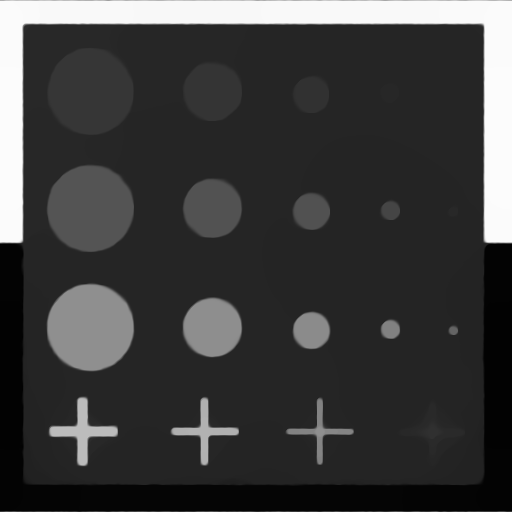}}};
	\spy[magnification=6] on (0.3125\textwidth,-0.306\textwidth) in node [name=c1]  at (0,-1.9\textwidth);
	\spy on (0.2825\textwidth,0.279\textwidth) in node [name=c1]  at (0,-0.95\textwidth);
	\spy on (-0.2925\textwidth,0.279\textwidth) in node [name=c1]  at (0,0\textwidth);
	\end{scope}
	\end{tikzpicture}
	\caption{\scalebox{1.2}{TV}}
	\label{fig:ImageDenoisingDeblurring-a}
	\end{subfigure}
	\centering
	    \begin{subfigure}[b]{0.18\textwidth}
	\centering \begin{tikzpicture}
	\begin{scope}[spy using outlines={rectangle,red,magnification=4,size=0.9\textwidth}]
	\node [name=c]{{\includegraphics[height=0.9\textwidth]{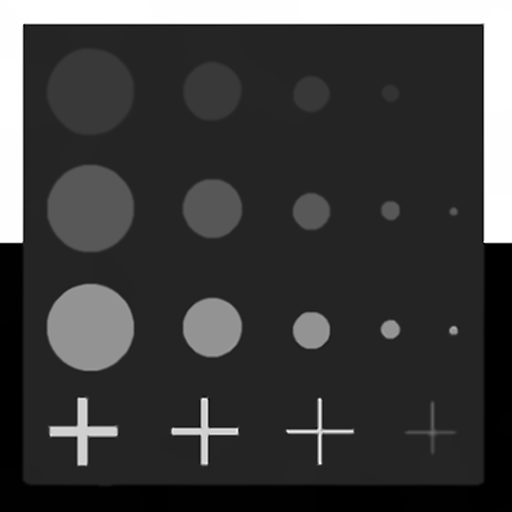}}};
	\spy[magnification=6] on (0.3125\textwidth,-0.306\textwidth) in node [name=c1]  at (0,-1.9\textwidth);
	\spy on (0.2825\textwidth,0.279\textwidth) in node [name=c1]  at (0,-0.95\textwidth);
	\spy on (-0.2925\textwidth,0.279\textwidth) in node [name=c1]  at (0,0\textwidth);
	\end{scope}
	\end{tikzpicture}
	\caption{\scalebox{1.2}{NLM}}
	\label{fig:ImageDenoisingDeblurring-b}
	\end{subfigure}
		\centering
	    \begin{subfigure}[b]{0.18\textwidth}
	\centering \begin{tikzpicture}
	\begin{scope}[spy using outlines={rectangle,red,magnification=4,size=0.9\textwidth}]
	\node [name=c]{{\includegraphics[height=0.9\textwidth]{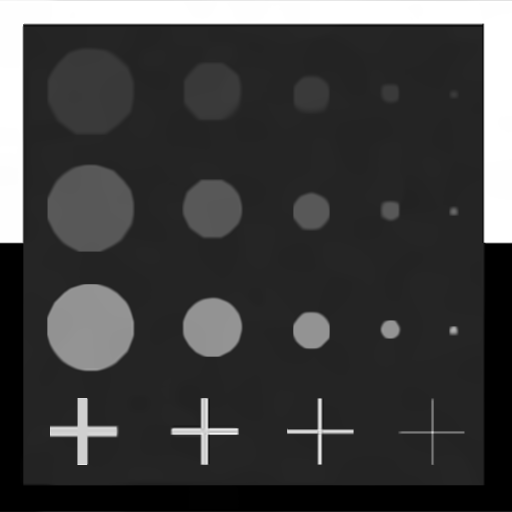}}};
	\spy[magnification=6] on (0.3125\textwidth,-0.306\textwidth) in node [name=c1]  at (0,-1.9\textwidth);
	\spy on (0.2825\textwidth,0.279\textwidth) in node [name=c1]  at (0,-0.95\textwidth);
	\spy on (-0.2925\textwidth,0.279\textwidth) in node [name=c1]  at (0,0\textwidth);
	\end{scope}
	\end{tikzpicture}
	\caption{\scalebox{1.2}{BM3D}}
	\label{fig:ImageDenoisingDeblurring-c}
	\end{subfigure}
	\centering
	\begin{subfigure}[b]{0.18\textwidth}
	\centering \begin{tikzpicture}
	\begin{scope}[spy using outlines={rectangle,red,magnification=4,size=0.9\textwidth}]
	\node [name=c]{{\includegraphics[height=0.9\textwidth]{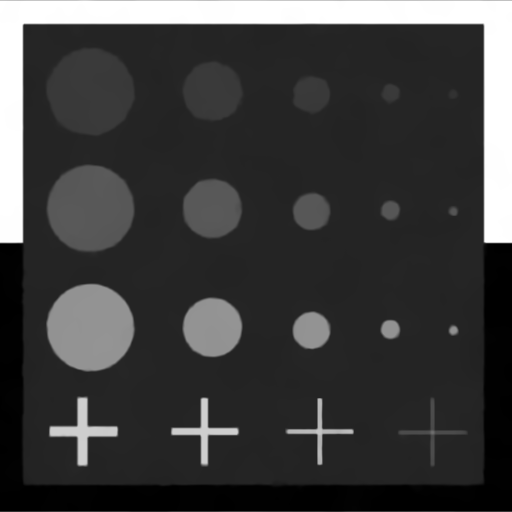}}};
	\spy[magnification=6] on (0.3125\textwidth,-0.306\textwidth) in node [name=c1]  at (0,-1.9\textwidth);
	\spy on (0.2825\textwidth,0.279\textwidth) in node [name=c1]  at (0,-0.95\textwidth);
	\spy on (-0.2925\textwidth,0.279\textwidth) in node [name=c1]  at (0,0\textwidth);
	\end{scope}
	\end{tikzpicture}
	\caption{\scalebox{1.2}{BM3D-WL1}}
	\label{fig:ImageDenoisingDeblurring-d}
	\end{subfigure}
	\centering
	
    \begin{subfigure}[b]{0.18\textwidth}
	\centering \begin{tikzpicture}
	\begin{scope}[spy using outlines={rectangle,red,magnification=4,size=0.9\textwidth}]
	\node [name=c]{{\includegraphics[height=0.9\textwidth]{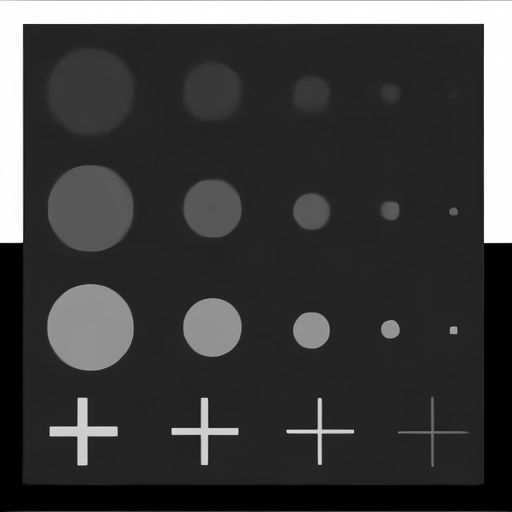}}};
	\spy[magnification=6] on (0.3125\textwidth,-0.306\textwidth) in node [name=c1]  at (0,-1.9\textwidth);
	\spy on (0.2825\textwidth,0.279\textwidth) in node [name=c1]  at (0,-0.95\textwidth);
	\spy on (-0.2925\textwidth,0.279\textwidth) in node [name=c1]  at (0,0\textwidth);
	\end{scope}
	\end{tikzpicture}
	\caption{\scalebox{1.2}{ICNN}}
	\label{fig:ImageDenoisingDeblurring-e}
	\end{subfigure}
    \begin{subfigure}[b]{0.18\textwidth}
	\centering \begin{tikzpicture}
	\begin{scope}[spy using outlines={rectangle,red,magnification=4,size=0.9\textwidth}]
	\node [name=c]{{\includegraphics[height=0.9\textwidth]{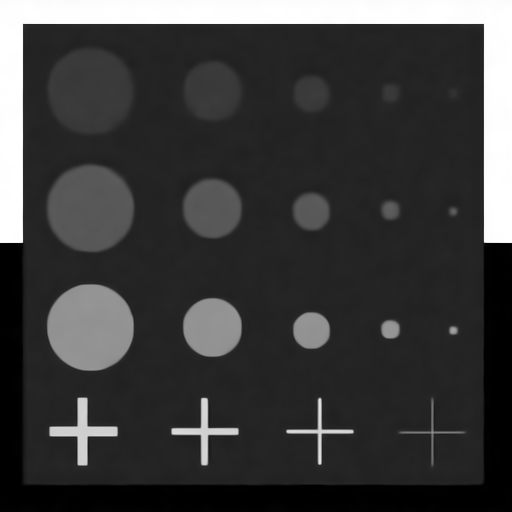}}};
	\spy[magnification=6] on (0.3125\textwidth,-0.306\textwidth) in node [name=c1]  at (0,-1.9\textwidth);
	\spy on (0.2825\textwidth,0.279\textwidth) in node [name=c1]  at (0,-0.95\textwidth);
	\spy on (-0.2925\textwidth,0.279\textwidth) in node [name=c1]  at (0,0\textwidth);
	\end{scope}
	\end{tikzpicture}
	\caption{\scalebox{1.2}{GCNN}}
	\label{fig:ImageDenoisingDeblurring-f}
	\end{subfigure}
    \begin{subfigure}[b]{0.18\textwidth}
	\centering \begin{tikzpicture}
	\begin{scope}[spy using outlines={rectangle,red,magnification=4,size=0.9\textwidth}]
	\node [name=c]{{\includegraphics[height=0.9\textwidth]{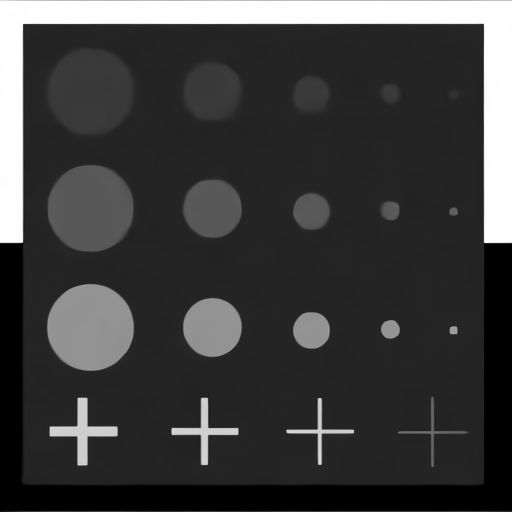}}};
	\spy[magnification=6] on (0.3125\textwidth,-0.306\textwidth) in node [name=c1]  at (0,-1.9\textwidth);
	\spy on (0.2825\textwidth,0.279\textwidth) in node [name=c1]  at (0,-0.95\textwidth);
	\spy on (-0.2925\textwidth,0.279\textwidth) in node [name=c1]  at (0,0\textwidth);
	\end{scope}
	\end{tikzpicture}
	\caption{\scalebox{1.2}{ICNN-TV}}
	\label{fig:ImageDenoisingDeblurring-g}
	\end{subfigure}
	\begin{subfigure}[b]{0.18\textwidth}
	\centering \begin{tikzpicture}
	\begin{scope}[spy using outlines={rectangle,red,magnification=4,size=0.9\textwidth}]
	\node [name=c]{{\includegraphics[height=0.9\textwidth]{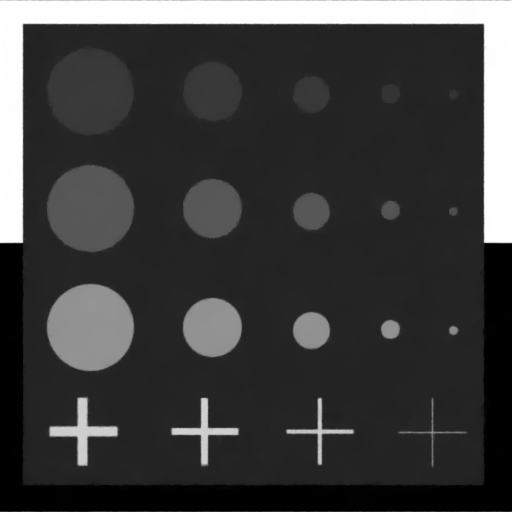}}};
	\spy[magnification=6] on (0.3125\textwidth,-0.306\textwidth) in node [name=c1]  at (0,-1.9\textwidth);
	\spy on (0.2825\textwidth,0.279\textwidth) in node [name=c1]  at (0,-0.95\textwidth);
	\spy on (-0.2925\textwidth,0.279\textwidth) in node [name=c1]  at (0,0\textwidth);
	\end{scope}
	\end{tikzpicture}
	\caption{\scalebox{1.2}{GCNN-TV}}
	\label{fig:ImageDenoisingDeblurring-h}
	\end{subfigure}}
	\caption{Three close-ups for each reconstruction by different methods obtained for the synthetic image.}
	\label{fig:ImageDenoisingDeblurring}
\end{figure}

In Figure \ref{fig:ImageDenoisingDeblurring}, for each method we report  the three restored zooms in the same range of gray levels. For what concerns the considered low-contrast circles, reported in the first two rows, it is evident that the hybrid approaches (such as BM3D-WL1,  ICNN-TV and GCNN-TV) outperform the other algorithms which exploit only one prior (TV, NLM, BM3D,ICNN, GCNN). Indeed, TV and NLM struggle to retrieve the small details pointed by the magenta arrows (Figures \ref{fig:ImageDenoisingDeblurring-a} and \ref{fig:ImageDenoisingDeblurring-b}), BM3D deforms the shape of the objects (Figure \ref{fig:ImageDenoisingDeblurring-c}), whereas  ICNN and GCNN do not provide solutions with sharp edges (Figures \ref{fig:ImageDenoisingDeblurring-e} and \ref{fig:ImageDenoisingDeblurring-f}). Focusing on the restoration of an object, the one-pixel thick cross, with a different shape and contrast,  we observe that BM3D, ICNN and GCNN achieve  the highest enhancement (see the last row of Figure \ref{fig:ImageDenoisingDeblurring}). 
However we remark  that, even in this case, TV and NLM tend to suppress very thin details.

\begin{figure}
\centering
    \begin{subfigure}[b]{0.45\textwidth}	%
	    \includegraphics[height=2cm]{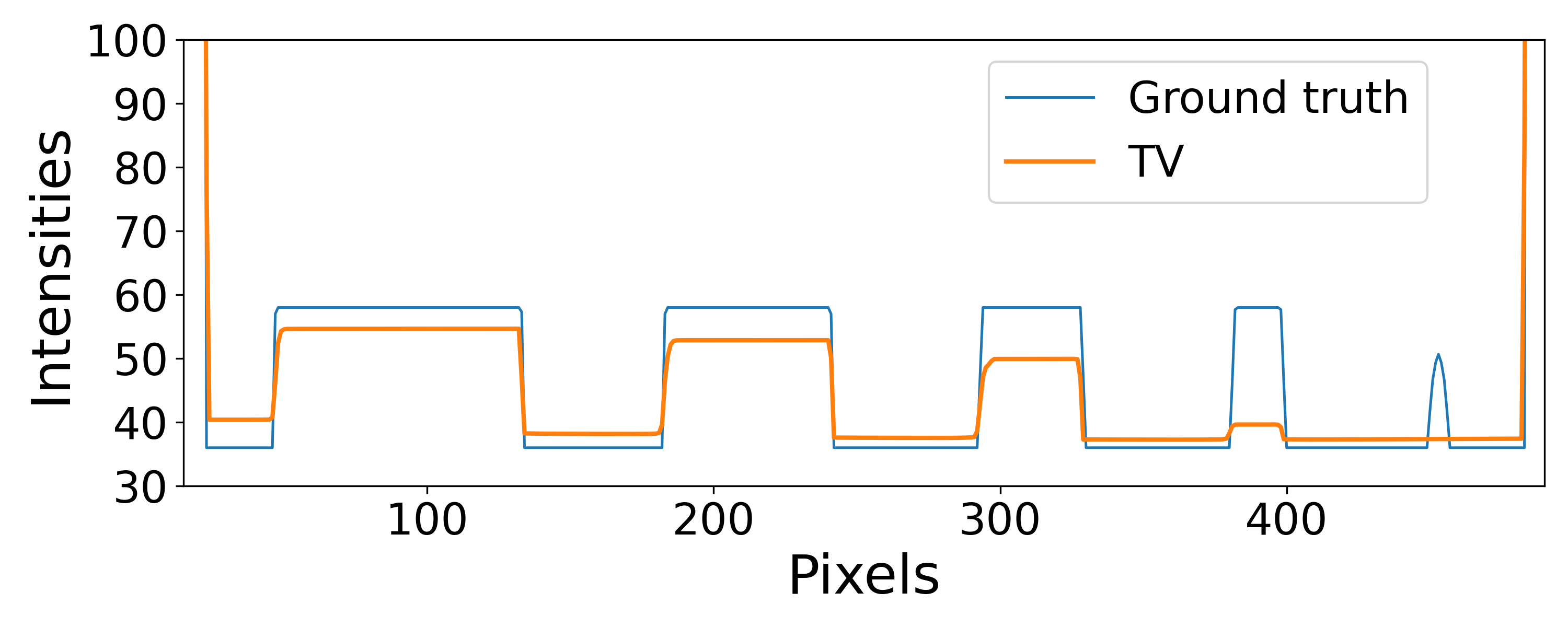}
	    \caption{TV}
	    \label{fig:ImageDenoisingDeblurring_profili1-a}
    \end{subfigure} \\ 
    \begin{subfigure}[b]{0.45\textwidth}
		\includegraphics[height=2cm]{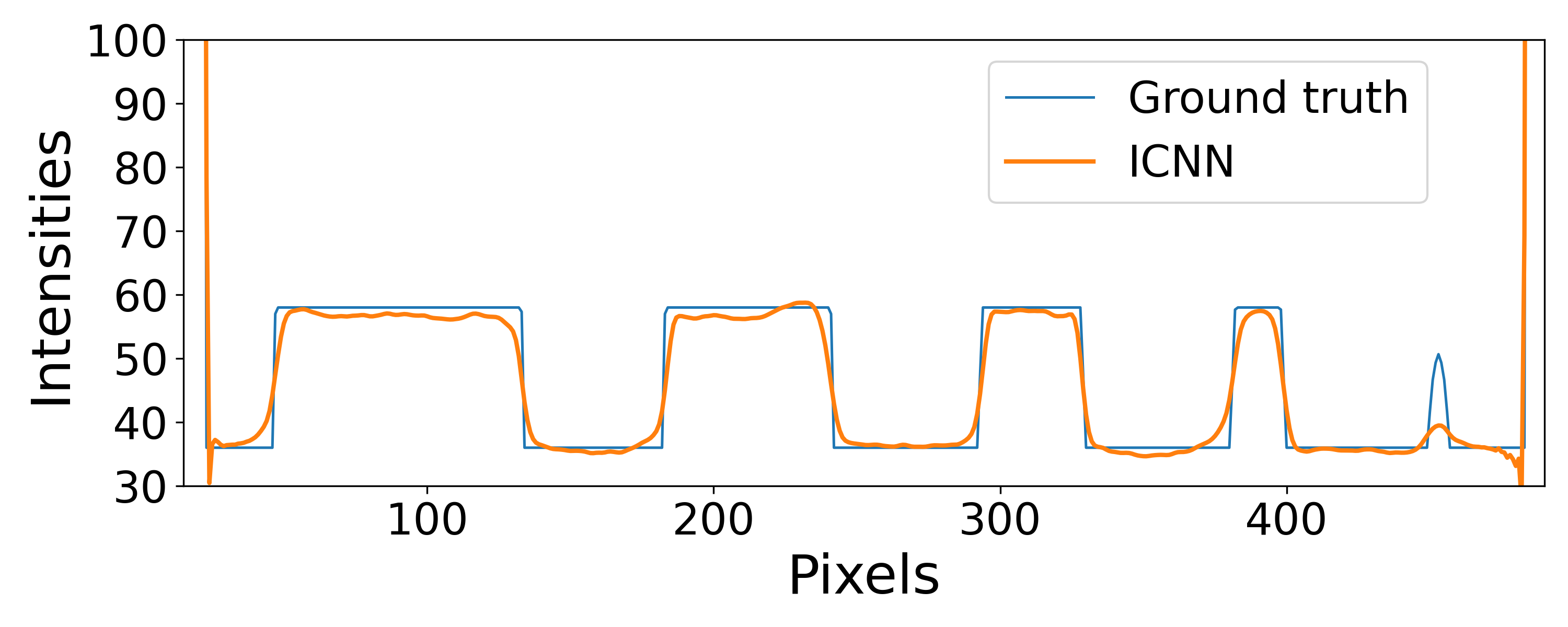}
		\caption{ICNN}
		\label{fig:ImageDenoisingDeblurring_profili1-b}
    \end{subfigure} \begin{subfigure}[b]{0.45\textwidth}
 		\includegraphics[height=2cm]{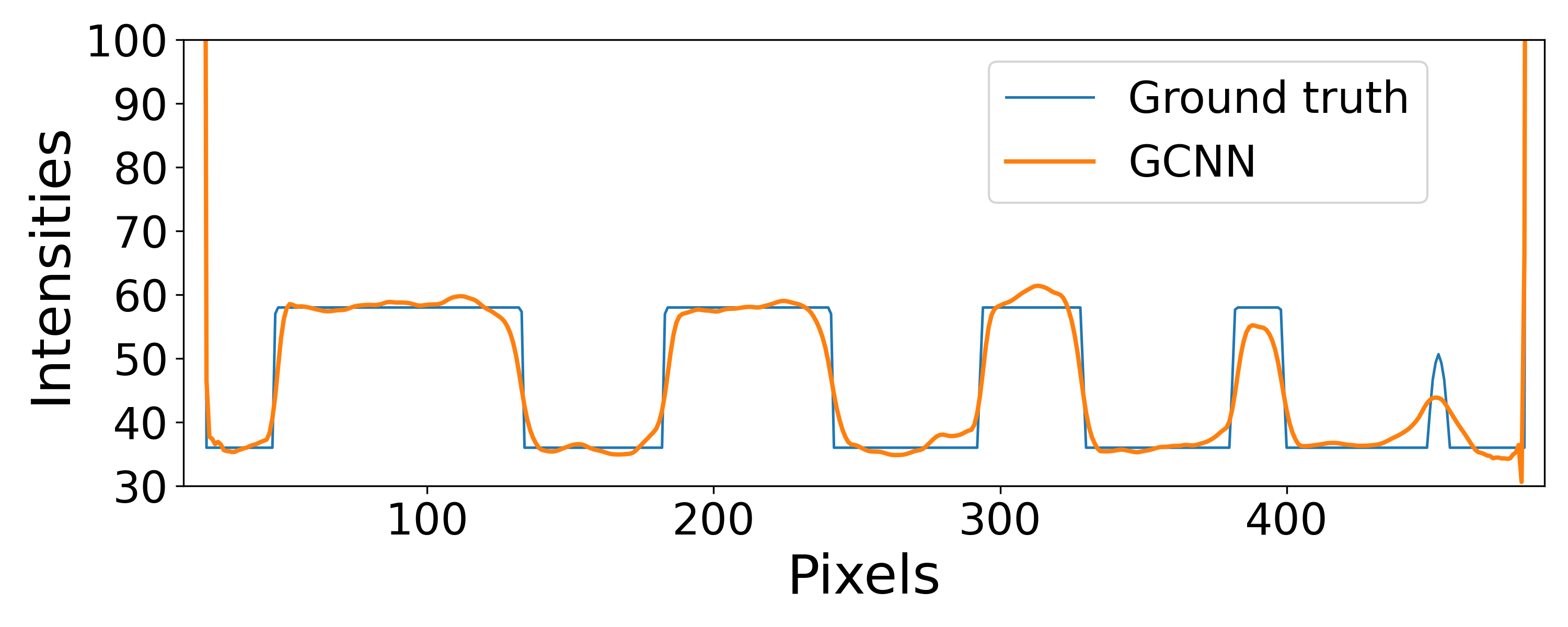}
 	\caption{GCNN}
 	\label{fig:ImageDenoisingDeblurring_profili1-c}
    \end{subfigure}\\
    \begin{subfigure}[b]{0.45\textwidth}
		\includegraphics[height=2cm]{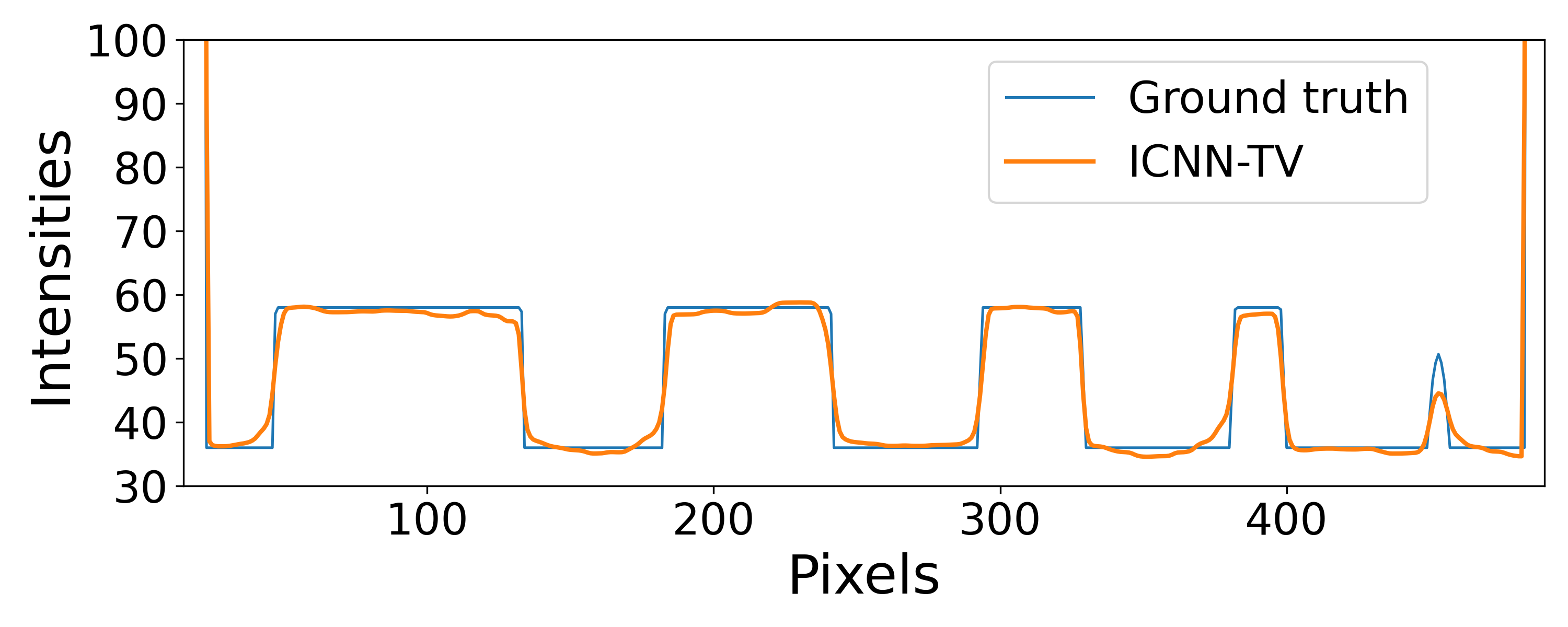}
	    \caption{ICNN-TV}
	    \label{fig:ImageDenoisingDeblurring_profili1-d}
    \end{subfigure}\begin{subfigure}[b]{0.45\textwidth}
 		\includegraphics[height=2cm]{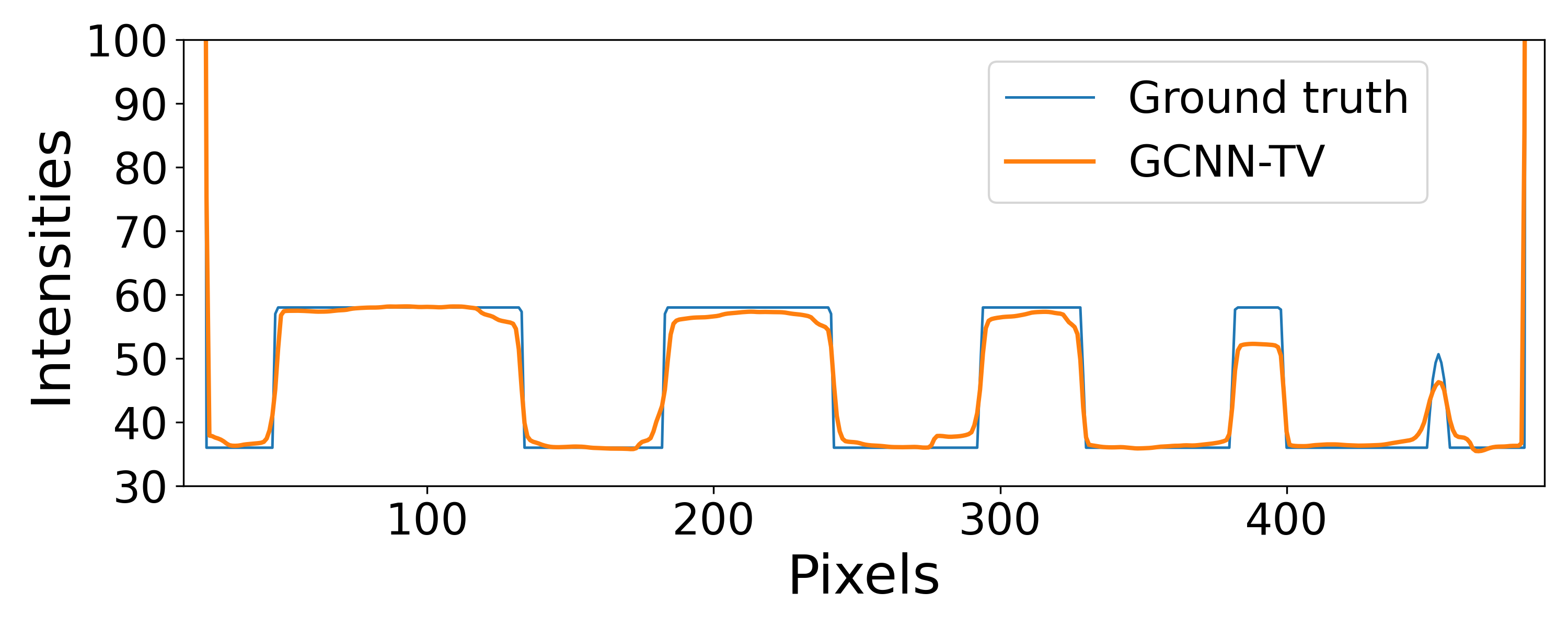}
 		\caption{GCNN-TV}
 		\label{fig:ImageDenoisingDeblurring_profili1-e}
    \end{subfigure}
	\caption{Intensity line profiles on the 90th row cutting the lowest contrasted circles. The blue and orange lines represent the ground truth and the restored image profiles for different methods, respectively.}
	\label{fig:ImageDenoisingDeblurring_profili1}
    \end{figure}

In Figure \ref{fig:ImageDenoisingDeblurring_profili1} we plot the pixel intensities of a  horizontal image row passing through all the lowest-contrasted circles, to  better inspect the effects of adding the TV internal prior to the ICNN and GCNN schemes on the most challenging objects.
The plot in Figure \ref{fig:ImageDenoisingDeblurring_profili1-a} reflects the typical loss-of-contrast drawback of the TV prior, oversmoothing the two smallest circles.
Adding the TV prior to ICNN and GCNN algorithms removes the residual noise, especially visible in the largest circle, while enhancing the edges. However, the  restoration of the smallest object has the major gain compared to the hybrid framework.

\begin{table}
\centering
\begin{adjustbox}{max width=\textwidth}
\scalebox{0.7}{
\begin{tabular}{l| ll | ll | ll }
\thickhline
{} &  \multicolumn{2}{c|}{AWGN of $std=10$} & \multicolumn{2}{c|}{AWGN of $std=15$} & \multicolumn{2}{c}{AWGN of $std=20$} \\
{} &    PSNR &    ROI-std &     ROI-PSNR &    ROI-std &     PSNR &    ROI-std  \\
\thickhline %
TV          &   30.8085         &  \TextGB{0.0271}  &  28.6664          & \TextGB{0.0507}   & 27.4028          & \TextGB{0.0772} \\
NLM     &    32.6266         &  \TextG{0.0896}   &   31.3772	& \TextG{0.1122}   &   30.1042		& \TextGB{0.1382} \\
BM3D     &    32.1221         & 0.3657    &   31.3283         &   0.5785  &    30.4806      &  0.7281\\
BM3D-WL1     &    31.7616         & 0.3974    &   30.7724         &   0.5802  &    30.2951      &  0.7779\\
ICNN       &   \TextB{34.1838}  &  0.4398           &  \TextB{33.0519}  & 0.5555            & \TextB{32.9788}  & 0.7851 \\
GCNN       &  \TextBB{34.7078} &  0.4749           &  \TextBB{33.9640} & 0.6568            & \TextBB{33.2446} & 0.8189 \\
ICNN-TV     &   32.3531         &  0.4081          &  31.3775          & 0.4798            & 30.4499          & 0.5553 \\
GCNN-TV     &   33.2648          &  0.1706   &  31.7743          & 0.2512    & 30.6453         & 0.3129 \\
\thickhline %
\end{tabular}}
\end{adjustbox}
\caption{Measures computed on restored images varying the standard deviation of the AWGN. The two best PSNR and ROI-std (standard deviation computed inside the green square in Figure \ref{fig:cirs-A})  values for each AWGN are highlighted in blue and green, respectively. The first best is highlighted in bold.}
\label{tab:Deblur1}
\end{table}

To test the robustness of the proposed models with respect to the noise, we analyze the results, reported in Table \ref{tab:Deblur1}, obtained by the considered methods when different variances of the AWGN are considered. We observe that, in terms of PSNR, the GCNN method gets the best values in all the cases, thus confirming the effectiveness, on this image, of the proposed CNN denoiser defined on the image gradient domain.
When we introduce the contribution of the TV-based internal prior, the PSNR values decrease, even if the global denoising effect due to TV is visually evident, as previously underlined. To confirm this, we report in Table \ref{tab:Deblur1} the standard deviation (ROI-std) computed on the constant region marked by the green bounding square in Figure \ref{fig:cirs-A}. The TV and NLM methods always have the lowest values, whereas the proposed hybrid approaches ICNN-TV and GCNN-TV are more effective in case of high noise. 
We can finally state that the combination of an external CNN denoiser and the internal TV prior  efficiently deals with noise removal and shape recovery.

\subsection{Results on real CT medical images} \label{subsec:ImageDenoisingDeblurring_CT}
We now consider two X-ray Computed Tomography  images to compare the effectiveness of the proposed schemes. In order to illustrate the advantages of our proposals, according to their features highlighted in the synthetic case, we examine a head and a chest CT image containing small and low-contrasted details.

\subsubsection{CT head image for epidural hemorrhage detection}

We now analyse a head tomographic image, downloaded from an open source dataset\footnote{\url{https://www.kaggle.com/vbookshelf/computed-tomography-ct-images}}, showing an intracranial hemorrhage and depicted in Figure \ref{fig:emo_a}. Intracranial hemorrhage, e.g. bleeding that occurs inside the cranium, requires often a rapid and intensive medical treatment based on the accurate localization of the  blood in the CT image obtained by segmentation algorithms. 
If the image is severely corrupted, the segmentation procedure may fail. As an example, we  blur the ground truth image with a Gaussian kernel of size $15 \times 15$ and standard deviation $0.5$ and we add AWGN with standard deviation $25$, thus we compute the  segmentation masks by an online open source software \footnote{\label{footnote:software} \url{http://brain.test.woza.work/}} shown in red in Figures \ref{fig:emo_a} and \ref{fig:emo_b}.
To highlight the importance of deblurring and denoising the image before segmenting it, we show the red mask computed on one restored image in Figure \ref{fig:emo_c}.

In Figure \ref{fig:emo_closeup} we report three close-ups for each method. The first one highlights the central part of the head CT image containing blood vessels, whereas the second zooms  shows a portion of the cerebral cortex with sulci. By a visual comparison, we observe that in both cases  TV, NLM, BM3D-WL1 output images look too smooth and blocky whereas the BM3D deforms the anatomical contours.
Focusing on the GCNN-based reconstructions, we highlight that the GCNN method accurately restores the vessels and sulci borders but it adds small artifacts which are significantly smoothed out in the GCNN-TV image. 
The third zoom of Figure \ref{fig:emo_closeup} focuses onto the epidural hemorrhage (pointed by the magenta arrow). We also report the widely used Jaccard similarity coefficient (J)  between the masks computed on the GT and the restored images: we highlight that the highest values are achieved by the three hybrid frameworks with an internal prior, i.e. BM3D-WL1, ICNN-TV and GCNN-TV.

\begin{figure}
    \centering
    \begin{subfigure}[b]{.3\textwidth}
    \centering
    \includegraphics[width=0.9\textwidth]{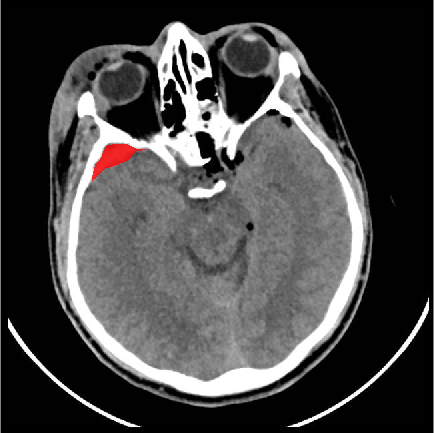}
    \caption{GT + mask}
    \label{fig:emo_a}
    \end{subfigure}\begin{subfigure}[b]{.3\textwidth}
    \centering
    \includegraphics[width=0.9\textwidth]{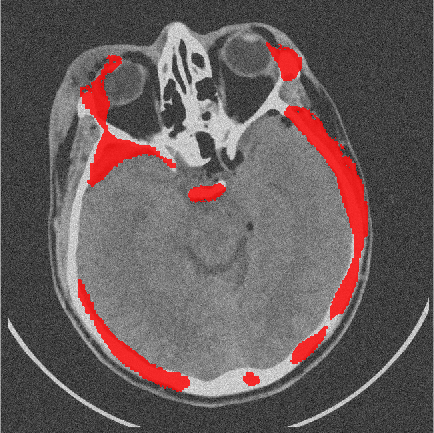}
    \caption{Corrupted + mask}
    \label{fig:emo_b}
    \end{subfigure} 
    \begin{subfigure}[b]{.3\textwidth}
    \includegraphics[width=0.9\textwidth]{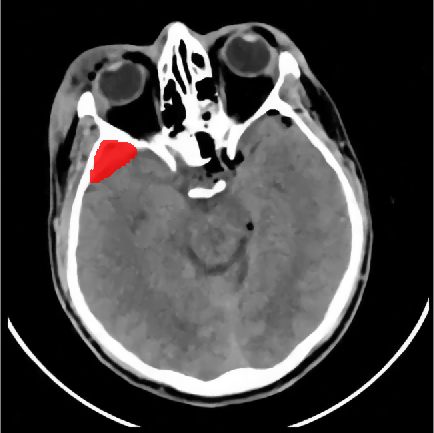}
    \caption{ICNN-TV + mask}
    \label{fig:emo_c}
    \end{subfigure}
    \caption{Head tomographic image with epidural hemorrhage. Computed masks are coloured red.}
    \label{fig:emo}
\end{figure}
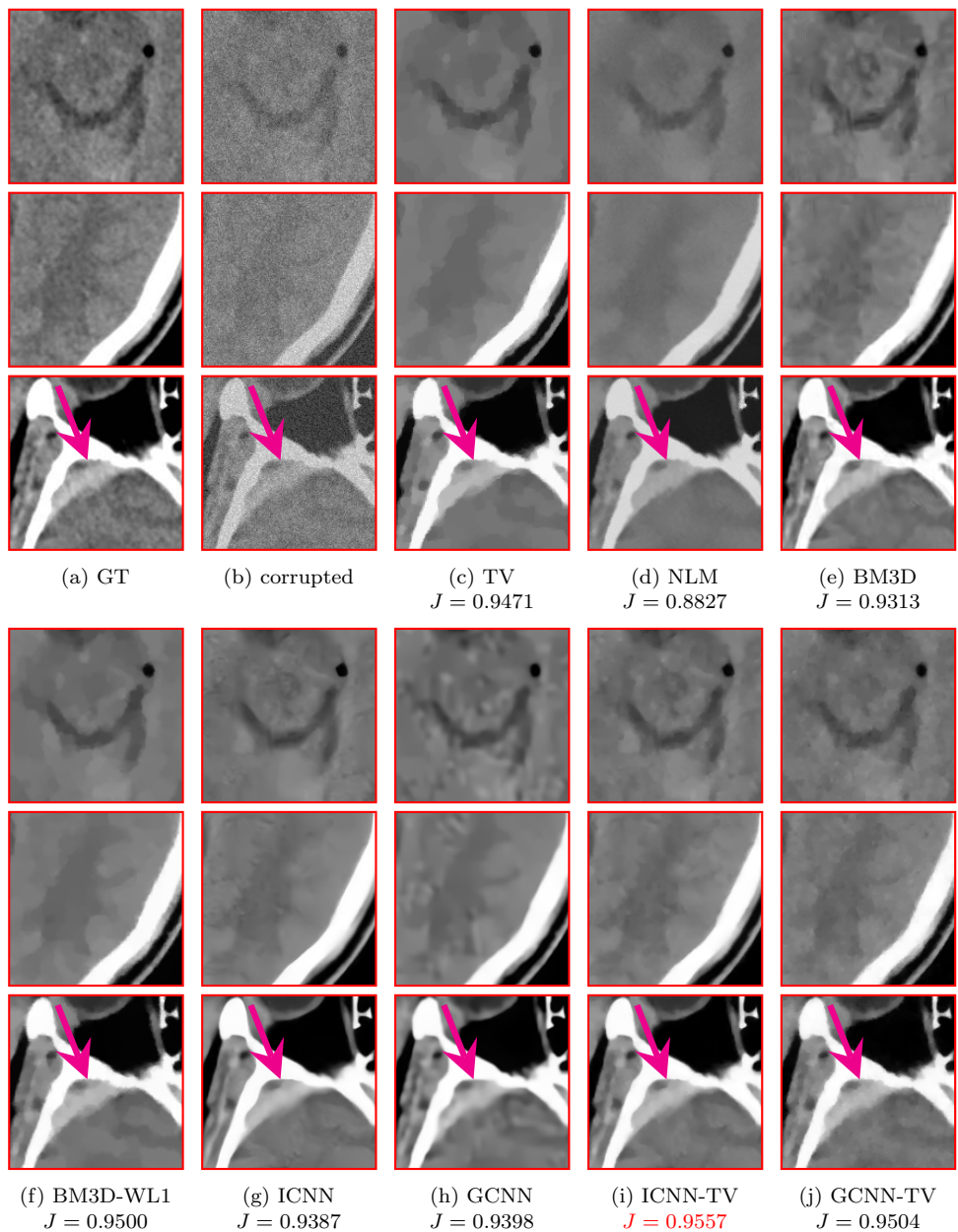
\begin{figure*}
\centering
	\begin{subfigure}[t]{0.19\textwidth}
		\captionsetup{justification=centering}
    \begin{tikzpicture}
	\begin{scope}[spy using outlines={rectangle,red,magnification=4,size=0.9\textwidth}]
	\node [name=c]{{\includegraphics[height=0.9\textwidth]{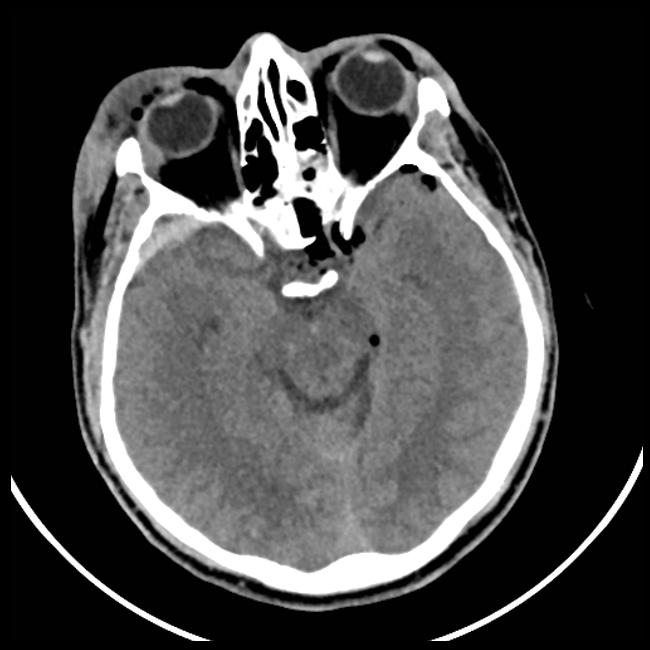}}};
	\spy [magnification=4] on (0.17\textwidth,-0.18\textwidth) in node [name=c1]  at (0,-0.95\textwidth);
    \spy[magnification=4] on (-0.2\textwidth,0.15\textwidth) in node [name=c1]  at (0,-1.90\textwidth);
    \spy [magnification=4] on (0\textwidth,-0.08\textwidth) in node [name=c1]  at (0,0\textwidth);
	\draw [-stealth, line width=0.6pt, magenta] (-0.65,0.65) -- ++(0.1,-0.25);
	\end{scope}
	\end{tikzpicture}
	\caption{{GT}}
	\label{fig:emo-a}
	\end{subfigure}\begin{subfigure}[t]{0.19\textwidth}	\captionsetup{justification=centering}
    \begin{tikzpicture}
	\begin{scope}[spy using outlines={rectangle,red,magnification=4,size=0.9\textwidth}]
	\node [name=c]{{\includegraphics[height=0.9\textwidth]{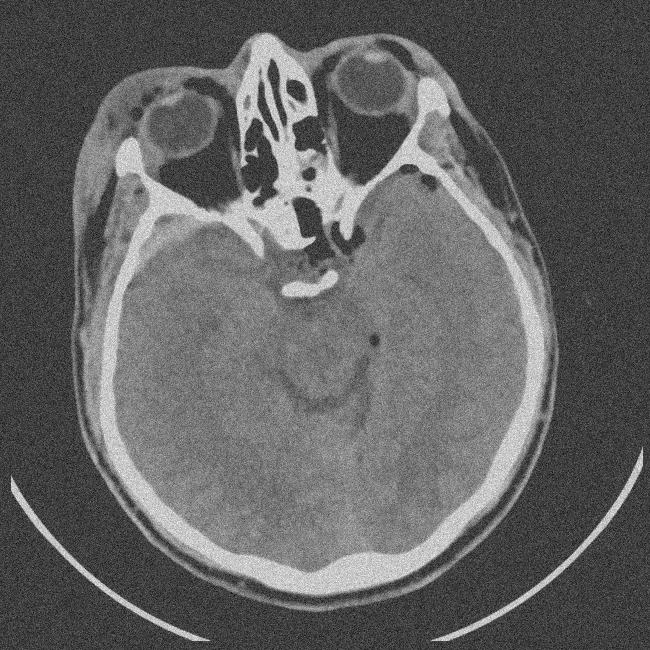}}};
   	\spy [magnification=4] on (0.17\textwidth,-0.18\textwidth) in node [name=c1]  at (0,-0.95\textwidth);
    \spy[magnification=4] on (-0.2\textwidth,0.15\textwidth) in node [name=c1]  at (0,-1.90\textwidth);
    \spy [magnification=4] on (0\textwidth,-0.08\textwidth) in node [name=c1]  at (0,0\textwidth);
	\draw [-stealth, line width=0.6pt, magenta] (-0.65,0.65) -- ++(0.1,-0.25);
	\end{scope}
	\end{tikzpicture}
	\caption{{corrupted}}
	\label{fig:emo-b}
	\end{subfigure}\begin{subfigure}[t]{0.19\textwidth}
		\captionsetup{justification=centering}
    \begin{tikzpicture}
	\begin{scope}[spy using outlines={rectangle,red,magnification=4,size=0.9\textwidth}]
	\node [name=c]{{\includegraphics[height=0.9\textwidth]{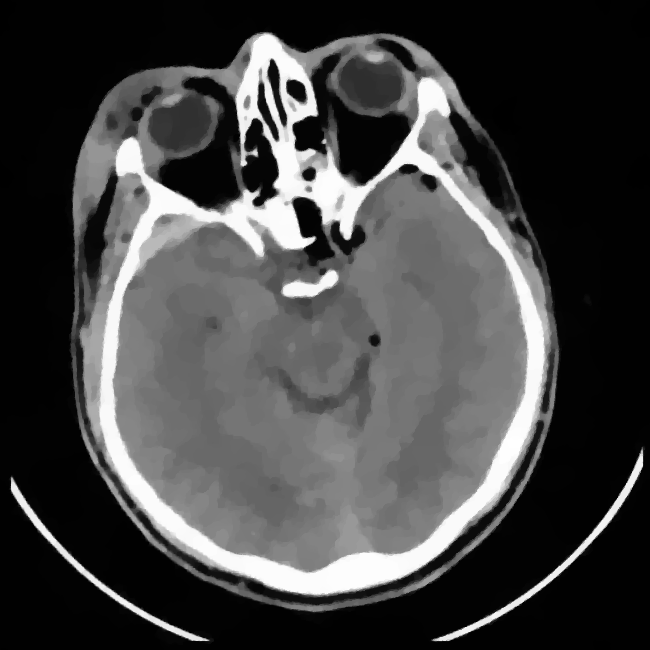}}};
   	\spy [magnification=4] on (0.17\textwidth,-0.18\textwidth) in node [name=c1]  at (0,-0.95\textwidth);
    \spy[magnification=4] on (-0.2\textwidth,0.15\textwidth) in node [name=c1]  at (0,-1.90\textwidth);
    \spy [magnification=4] on (0\textwidth,-0.08\textwidth) in node [name=c1]  at (0,0\textwidth);
	\draw [-stealth, line width=0.6pt, magenta] (-0.65,0.65) -- ++(0.1,-0.25);
	\end{scope}
	\end{tikzpicture}
	\caption{{TV \\ $J= 0.9471$}}
	\label{fig:emo-c}
	\end{subfigure}\begin{subfigure}[t]{0.19\textwidth}
		\captionsetup{justification=centering}
    \begin{tikzpicture}
	\begin{scope}[spy using outlines={rectangle,red,magnification=4,size=0.9\textwidth}]
	\node [name=c]{{\includegraphics[height=0.9\textwidth]{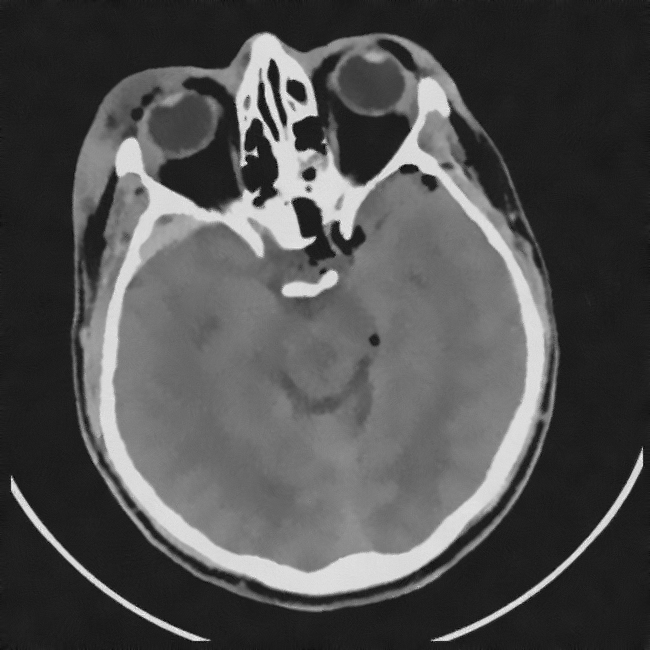}}};
    \spy [magnification=4] on (0.17\textwidth,-0.18\textwidth) in node [name=c1]  at (0,-0.95\textwidth);
    \spy[magnification=4] on (-0.2\textwidth,0.15\textwidth) in node [name=c1]  at (0,-1.90\textwidth);
    \spy [magnification=4] on (0\textwidth,-0.08\textwidth) in node [name=c1]  at (0,0\textwidth);
	\draw [-stealth, line width=0.6pt, magenta] (-0.65,0.65) -- ++(0.1,-0.25);
	\end{scope}
	\end{tikzpicture}
	\caption{{NLM \\ $J= 0.8827$}}
	\label{fig:emo-d}
	\end{subfigure}\begin{subfigure}[t]{0.19\textwidth}
	\captionsetup{justification=centering}
    \begin{tikzpicture}
	\begin{scope}[spy using outlines={rectangle,red,magnification=4,size=0.9\textwidth}]
	\node [name=c]{{\includegraphics[height=0.9\textwidth]{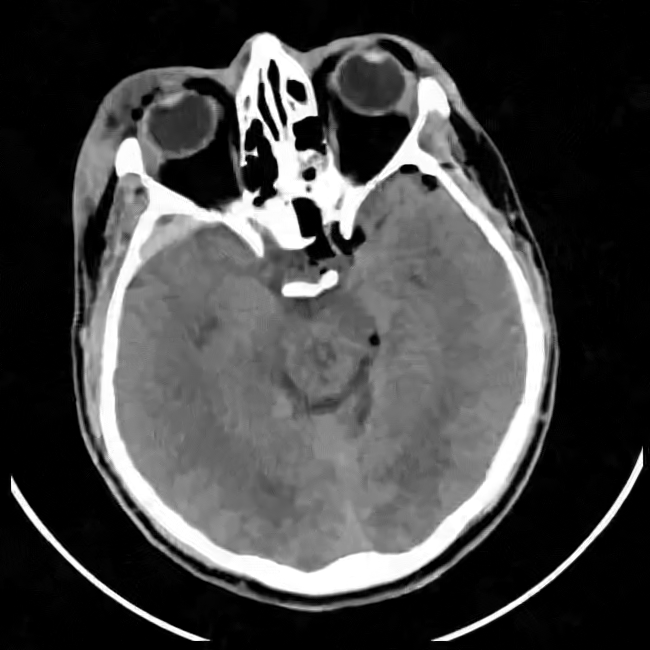}}};
   	\spy [magnification=4] on (0.17\textwidth,-0.18\textwidth) in node [name=c1]  at (0,-0.95\textwidth);
    \spy[magnification=4] on (-0.2\textwidth,0.15\textwidth) in node [name=c1]  at (0,-1.90\textwidth);
    \spy [magnification=4] on (0\textwidth,-0.08\textwidth) in node [name=c1]  at (0,0\textwidth);
	\draw [-stealth, line width=0.6pt, magenta] (-0.65,0.65) -- ++(0.1,-0.25);
	\end{scope}
	\end{tikzpicture}
	\caption{ BM3D\\$J=0.9313$} 
	\label{fig:emo-e}
	\end{subfigure}
	
	\begin{subfigure}[t]{0.19\textwidth}
	\captionsetup{justification=centering}
    \begin{tikzpicture}
	\begin{scope}[spy using outlines={rectangle,red,magnification=4,size=0.9\textwidth}]
	\node [name=c]{{\includegraphics[height=0.9\textwidth]{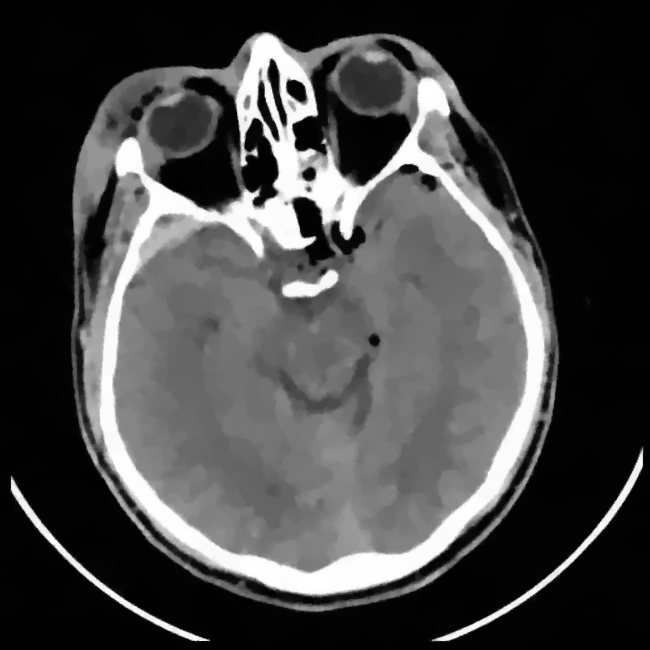}}};
     	\spy [magnification=4] on (0.17\textwidth,-0.18\textwidth) in node [name=c1]  at (0,-0.95\textwidth);
    \spy[magnification=4] on (-0.2\textwidth,0.15\textwidth) in node [name=c1]  at (0,-1.90\textwidth);
    \spy [magnification=4] on (0\textwidth,-0.08\textwidth) in node [name=c1]  at (0,0\textwidth);
	\draw [-stealth, line width=0.6pt, magenta] (-0.65,0.65) -- ++(0.1,-0.25);
	\end{scope}
	\end{tikzpicture}
	\caption{{BM3D-WL1 \\ $J= 0.9500$}} %
	\label{fig:emo-f}
	\end{subfigure}\begin{subfigure}[t]{0.19\textwidth}	\captionsetup{justification=centering}
	\begin{tikzpicture}
	\begin{scope}[spy using outlines={rectangle,red,magnification=4,size=0.9\textwidth}]
	\node [name=c]{{\includegraphics[height=0.9\textwidth]{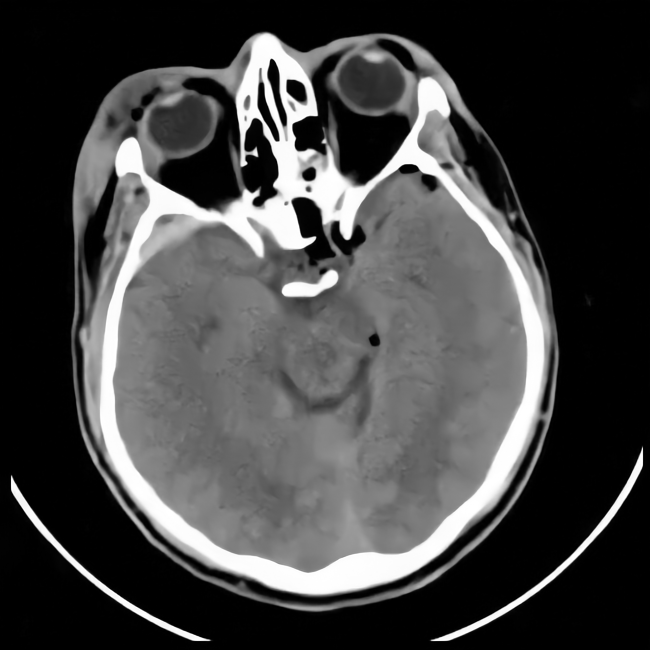}}};
    	\spy [magnification=4] on (0.17\textwidth,-0.18\textwidth) in node [name=c1]  at (0,-0.95\textwidth);
    \spy[magnification=4] on (-0.2\textwidth,0.15\textwidth) in node [name=c1]  at (0,-1.90\textwidth);
    \spy [magnification=4] on (0\textwidth,-0.08\textwidth) in node [name=c1]  at (0,0\textwidth);
	\draw [-stealth, line width=0.6pt, magenta] (-0.65,0.65) -- ++(0.1,-0.25);
	\end{scope}
	\end{tikzpicture}
	\caption{{ICNN \\ $J= 0.9387$}}
	\label{fig:emo-g}
	\end{subfigure}\begin{subfigure}[t]{0.19\textwidth}	\captionsetup{justification=centering}
	\begin{tikzpicture}
	\begin{scope}[spy using outlines={rectangle,red,magnification=4,size=0.9\textwidth}]
	\node [name=c]{{\includegraphics[height=0.9\textwidth]{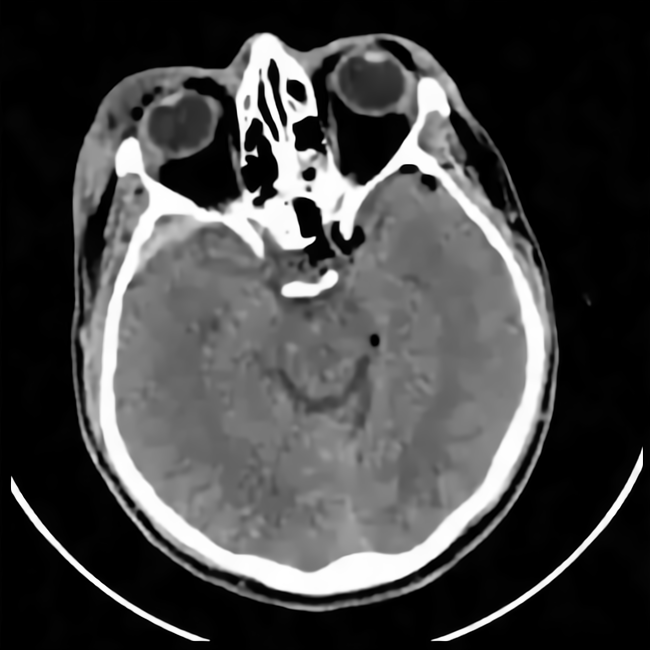}}};
    	\spy [magnification=4] on (0.17\textwidth,-0.18\textwidth) in node [name=c1]  at (0,-0.95\textwidth);
    \spy[magnification=4] on (-0.2\textwidth,0.15\textwidth) in node [name=c1]  at (0,-1.90\textwidth);
    \spy [magnification=4] on (0\textwidth,-0.08\textwidth) in node [name=c1]  at (0,0\textwidth);
	\draw [-stealth, line width=0.6pt, magenta] (-0.65,0.65) -- ++(0.1,-0.25);
	\end{scope}
	\end{tikzpicture}
	\caption{{GCNN \\ $J=0.9398$}}
	\label{fig:emo-h}
	\end{subfigure}\begin{subfigure}[t]{0.19\textwidth} 	\captionsetup{justification=centering}
    \begin{tikzpicture}
	\begin{scope}[spy using outlines={rectangle,red,magnification=4,size=0.9\textwidth}]
	\node [name=c]{{\includegraphics[height=0.9\textwidth]{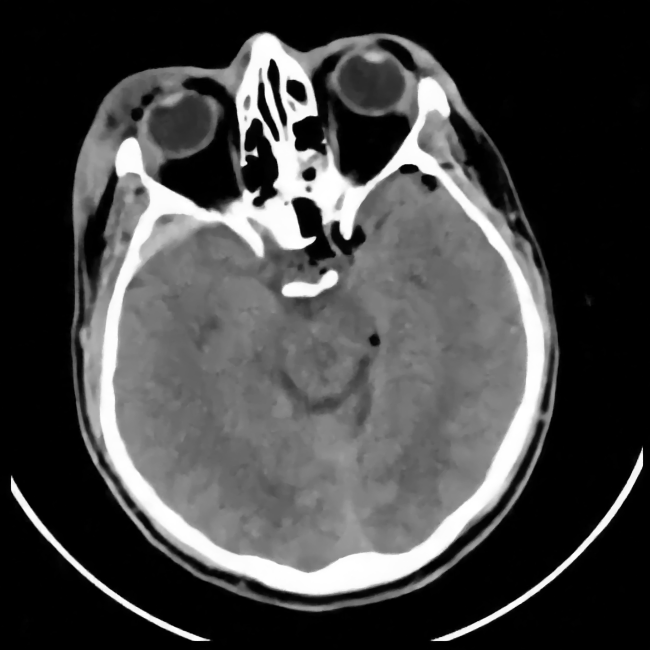}}};
	\spy [magnification=4] on (0.17\textwidth,-0.18\textwidth) in node [name=c1]  at (0,-0.95\textwidth);
    \spy[magnification=4] on (-0.2\textwidth,0.15\textwidth) in node [name=c1]  at (0,-1.90\textwidth);
    \spy [magnification=4] on (0\textwidth,-0.08\textwidth) in node [name=c1]  at (0,0\textwidth);
	\draw [-stealth, line width=0.6pt, magenta] (-0.65,0.65) -- ++(0.1,-0.25);
	\end{scope}
	\end{tikzpicture}
	\caption{{ICNN-TV \\ $\textcolor{red}{J=0.9557}$}}
	\label{fig:emo-i}
	\end{subfigure}\begin{subfigure}[t]{0.19\textwidth}	\captionsetup{justification=centering}
	\begin{tikzpicture}
	\begin{scope}[spy using outlines={rectangle,red,magnification=4,size=0.9\textwidth}]
	\node [name=c]{{\includegraphics[height=0.9\textwidth]{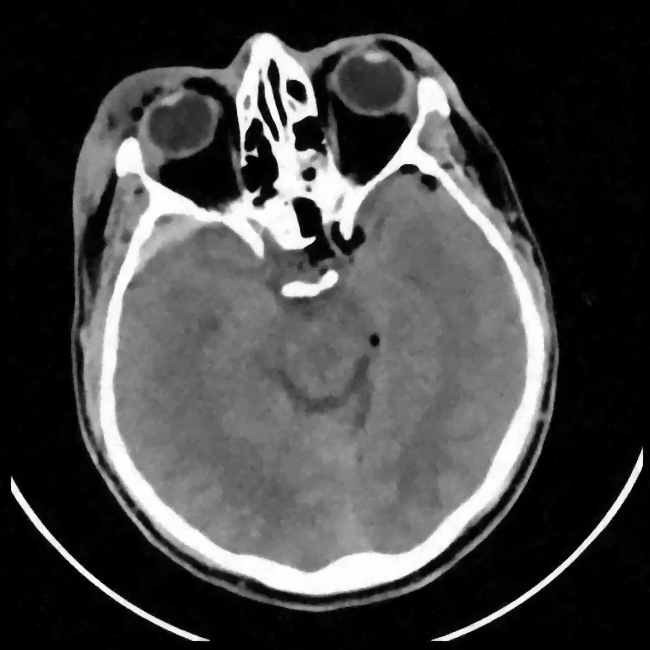}}};
   	\spy [magnification=4] on (0.17\textwidth,-0.18\textwidth) in node [name=c1]  at (0,-0.95\textwidth);
    \spy[magnification=4] on (-0.2\textwidth,0.15\textwidth) in node [name=c1]  at (0,-1.90\textwidth);
    \spy [magnification=4] on (0\textwidth,-0.08\textwidth) in node [name=c1]  at (0,0\textwidth);
	\draw [-stealth, line width=0.6pt, magenta] (-0.65,0.65) -- ++(0.1,-0.25);
	\end{scope}
	\end{tikzpicture}
	\caption{{GCNN-TV \\ $J=0.9504$}}
	\label{fig:emo-j}
	\end{subfigure}
	\caption{Three close-ups for each reconstruction by different methods obtained for the head CT image. The magenta arrows highlight the epidural hemorrhages. We highlight in red the best Jaccard index.}
	\label{fig:emo_closeup}
\end{figure*}

\subsubsection{Low-dose CT chest image}

\begin{figure}[ht]
\centering
    \begin{subfigure}[b]{.4\textwidth}
	\centering
	\begin{tikzpicture}
	\begin{scope}[spy using outlines={rectangle,red,magnification=2,size=1.5cm}]
	\node [name=c]{\frame{\includegraphics[width=0.95\textwidth]{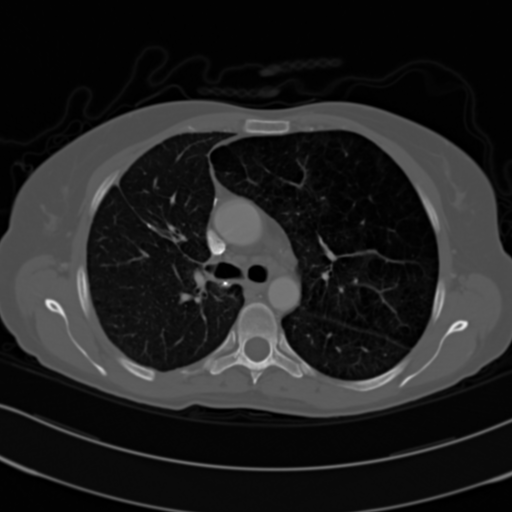}}};
	\draw[green,line width=0.1mm] (0.4\textwidth,0.4\textwidth) rectangle (0.3\textwidth,0.3\textwidth);
	\end{scope}
	\end{tikzpicture}
	\caption{Ground truth}
	\label{fig:low_dose_a}
	\end{subfigure}	\begin{subfigure}[b]{.4\textwidth}
	\centering
	\begin{tikzpicture}
	\begin{scope}[spy using outlines={rectangle,red,magnification=2,size=1.5cm}]
	\node [name=c]{\frame{\includegraphics[width=0.95\textwidth]{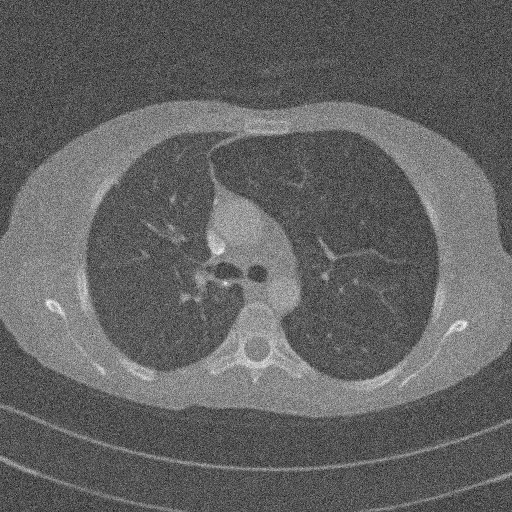}}};
	\end{scope}
	\end{tikzpicture}
	
	\caption{Corrupted}
	\label{fig:low_dose_b}
	\end{subfigure}
	
    \caption{Low-dose CT chest image (\texttt{ID: 0005}). In (a) the green square highlights the uniform patch used to evaluate ROI-std.}
     \label{fig:low_dose1}
    \end{figure}

We now consider a low-dose Computed Tomography chest image downloaded from an open source dataset\footnote{\url{https://www.kaggle.com/kmader/siim-medical-images}}. The considered image (ID: 0005) is depicted in Figure \ref{fig:low_dose_a}: the presence of soft tissues and organs yields uniform, flat and low-contrasted regions in the gray scale image. \\ %
We add blur by using a Gaussian kernel of dimension $15 \times 15$ with standard deviation $0.5$ and AWGN with standard deviation equals to 25. In Figure \ref{fig:low_dose_b} we show the corrupted image where small and low-contrasted details are not well detectable. 

\begin{figure*}
\centering
	\begin{subfigure}[b]{0.19\textwidth}
	\centering
    \begin{tikzpicture}
	\begin{scope}[spy using outlines={rectangle,red,magnification=4,size=0.9\textwidth}]
	\node [name=c]{{\includegraphics[height=0.9\textwidth]{results/low_dose2/low_dose2.png}}};
    \spy[magnification=5] on (0.184\textwidth,-0.02\textwidth) in node [name=c1]  at (0,-0.95\textwidth);
	\spy [magnification=4] on (0\textwidth,-0.18\textwidth) in node [name=c1]  at (0,-1.9\textwidth);
	 \spy [magnification=6] on (-0.05\textwidth,0.05\textwidth) in node [name=c1]  at (0,0\textwidth);
    \draw [-stealth, line width=0.1pt, magenta] (-0.02,-0.05) -- ++(-0.1,0.09);
	\end{scope}
	\end{tikzpicture}
	\caption{{GT}}
	\label{fig:low_dose-a}
	\end{subfigure}\begin{subfigure}[b]{0.19\textwidth}
	\centering
    \begin{tikzpicture}
	\begin{scope}[spy using outlines={rectangle,red,magnification=4,size=0.9\textwidth}]
	\node [name=c]{{\includegraphics[height=0.9\textwidth]{results/low_dose2/starting_low_dose2.png}}};
    \spy[magnification=5] on (0.184\textwidth,-0.02\textwidth) in node [name=c1]  at (0,-0.95\textwidth);
	\spy [magnification=4] on (0\textwidth,-0.18\textwidth) in node [name=c1]  at (0,-1.9\textwidth);
	 \spy [magnification=6] on (-0.05\textwidth,0.05\textwidth) in node [name=c1]  at (0,0\textwidth);
    \draw [-stealth, line width=0.1pt, magenta] (-0.02,-0.05) -- ++(-0.1,0.09);
	\end{scope}
	\end{tikzpicture}
	\caption{{corrupted}}
	\label{fig:low_dose-b}
	\end{subfigure}\begin{subfigure}[b]{0.19\textwidth}
	\centering
    \begin{tikzpicture}
	\begin{scope}[spy using outlines={rectangle,red,magnification=4,size=0.9\textwidth}]
	\node [name=c]{{\includegraphics[height=0.9\textwidth]{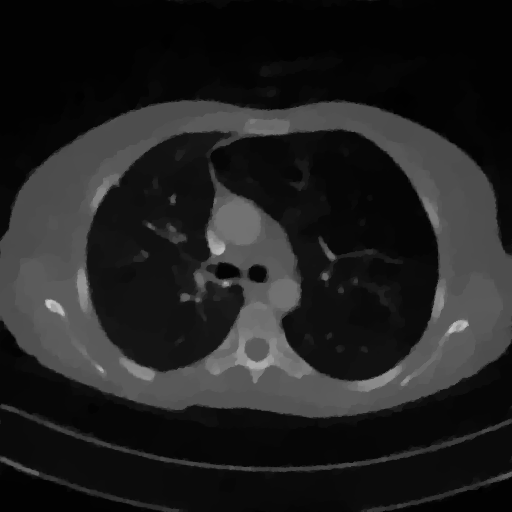}}};
     \spy[magnification=5] on (0.184\textwidth,-0.02\textwidth) in node [name=c1]  at (0,-0.95\textwidth);
	\spy [magnification=4] on (0\textwidth,-0.18\textwidth) in node [name=c1]  at (0,-1.9\textwidth);
	 \spy [magnification=6] on (-0.05\textwidth,0.05\textwidth) in node [name=c1]  at (0,0\textwidth);
    \draw [-stealth, line width=0.1pt, magenta] (-0.02,-0.05) -- ++(-0.1,0.09);
	\end{scope}
	\end{tikzpicture}
	\caption{{TV}}
	\label{fig:low_dose-c}
	\end{subfigure}\begin{subfigure}[b]{0.19\textwidth}
	\centering
    \begin{tikzpicture}
	\begin{scope}[spy using outlines={rectangle,red,magnification=4,size=0.9\textwidth}]
	\node [name=c]{{\includegraphics[height=0.9\textwidth]{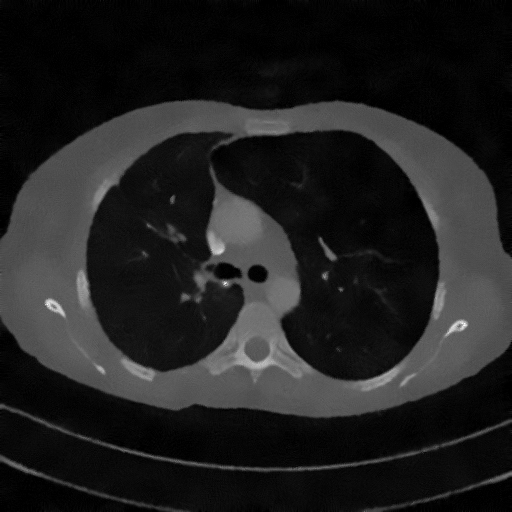}}};
    \spy[magnification=5] on (0.184\textwidth,-0.02\textwidth) in node [name=c1]  at (0,-0.95\textwidth);
	\spy [magnification=4] on (0\textwidth,-0.18\textwidth) in node [name=c1]  at (0,-1.9\textwidth);
	 \spy [magnification=6] on (-0.05\textwidth,0.05\textwidth) in node [name=c1]  at (0,0\textwidth);
    \draw [-stealth, line width=0.1pt, magenta] (-0.02,-0.05) -- ++(-0.1,0.09);
	\end{scope}
	\end{tikzpicture}
	\caption{{NLM}}
	\label{fig:low_dose-d}
	\end{subfigure}\begin{subfigure}[b]{0.19\textwidth}
	\centering 
	\begin{tikzpicture}
	\begin{scope}[spy using outlines={rectangle,red,magnification=4,size=0.9\textwidth}]
	\node [name=c]{{\includegraphics[height=0.9\textwidth]{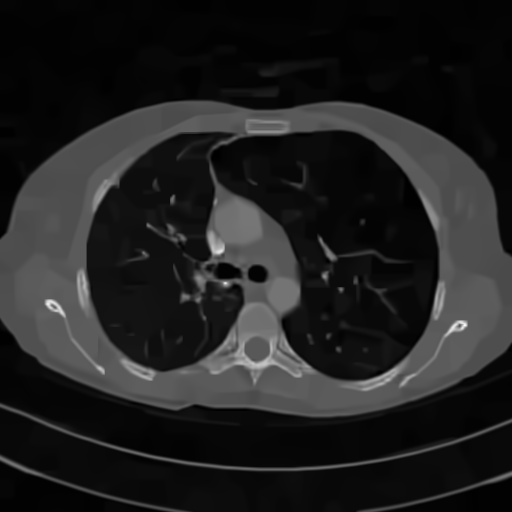}}};
    \spy[magnification=5] on (0.184\textwidth,-0.02\textwidth) in node [name=c1]  at (0,-0.95\textwidth);
	\spy [magnification=4] on (0\textwidth,-0.18\textwidth) in node [name=c1]  at (0,-1.9\textwidth);
	 \spy [magnification=6] on (-0.05\textwidth,0.05\textwidth) in node [name=c1]  at (0,0\textwidth);
    \draw [-stealth, line width=0.1pt, magenta] (-0.02,-0.05) -- ++(-0.1,0.09);
	\end{scope}
	\end{tikzpicture}
	\caption{{BM3D}}
	\label{fig:low_dose-e}
	\end{subfigure}
	\begin{subfigure}[b]{0.19\textwidth}
	\centering 
	\begin{tikzpicture}
	\begin{scope}[spy using outlines={rectangle,red,magnification=4,size=0.9\textwidth}]
	\node [name=c]{{\includegraphics[height=0.9\textwidth]{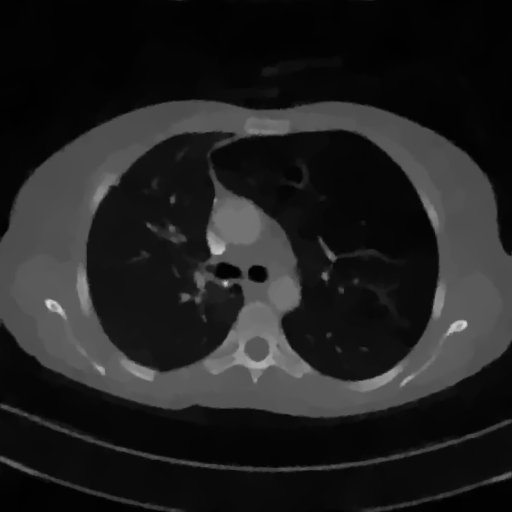}}};
    \spy[magnification=5] on (0.184\textwidth,-0.02\textwidth) in node [name=c1]  at (0,-0.95\textwidth);
	\spy [magnification=4] on (0\textwidth,-0.18\textwidth) in node [name=c1]  at (0,-1.9\textwidth);
	 \spy [magnification=6] on (-0.05\textwidth,0.05\textwidth) in node [name=c1]  at (0,0\textwidth);
    \draw [-stealth, line width=0.1pt, magenta] (-0.02,-0.05) -- ++(-0.1,0.09);
	\end{scope}
	\end{tikzpicture}
	\caption{{BM3D-WL1}}
	\label{fig:low_dose-f}
	\end{subfigure}\begin{subfigure}[b]{0.19\textwidth}
	\centering \begin{tikzpicture}
	\begin{scope}[spy using outlines={rectangle,red,magnification=4,size=0.9\textwidth}]
	\node [name=c]{{\includegraphics[height=0.9\textwidth]{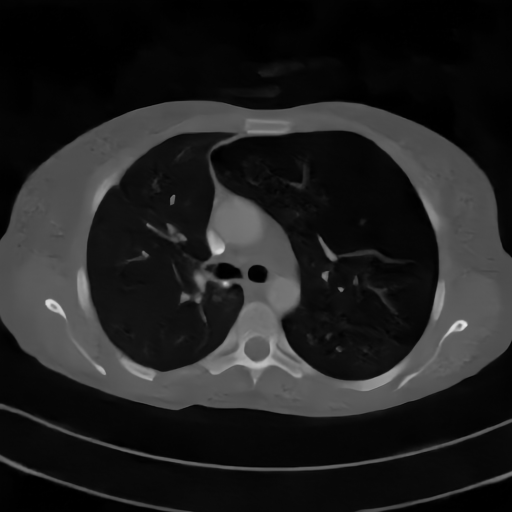}}};
    \spy[magnification=5] on (0.184\textwidth,-0.02\textwidth) in node [name=c1]  at (0,-0.95\textwidth);
	\spy [magnification=4] on (0\textwidth,-0.18\textwidth) in node [name=c1]  at (0,-1.9\textwidth);
	 \spy [magnification=6] on (-0.05\textwidth,0.05\textwidth) in node [name=c1]  at (0,0\textwidth);
    \draw [-stealth, line width=0.1pt, magenta] (-0.02,-0.05) -- ++(-0.1,0.09);
	\end{scope}
	\end{tikzpicture}
	\caption{{ICNN}}
	\label{fig:low_dose-g}
	\end{subfigure}\begin{subfigure}[b]{0.19\textwidth}
	\centering 
	\begin{tikzpicture}
	\begin{scope}[spy using outlines={rectangle,red,magnification=4,size=0.9\textwidth}]
	\node [name=c]{{\includegraphics[height=0.9\textwidth]{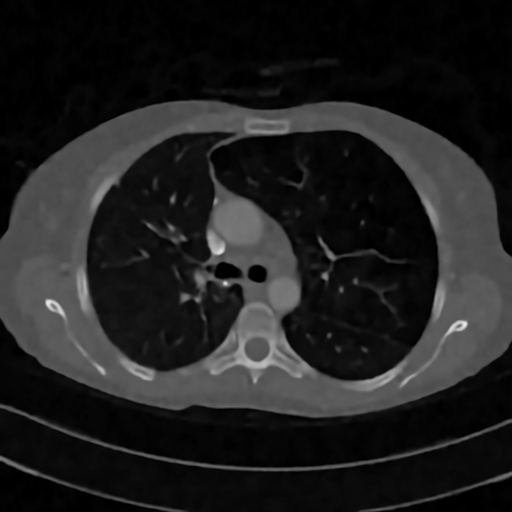}}};
     \spy[magnification=5] on (0.184\textwidth,-0.02\textwidth) in node [name=c1]  at (0,-0.95\textwidth);
	\spy [magnification=4] on (0\textwidth,-0.18\textwidth) in node [name=c1]  at (0,-1.9\textwidth);
	 \spy [magnification=6] on (-0.05\textwidth,0.05\textwidth) in node [name=c1]  at (0,0\textwidth);
    \draw [-stealth, line width=0.1pt, magenta] (-0.02,-0.05) -- ++(-0.1,0.09);
	\end{scope}
	\end{tikzpicture}
	\caption{{GCNN}}
	\label{fig:low_dose-h}
	\end{subfigure}\begin{subfigure}[b]{0.19\textwidth}
	\centering \begin{tikzpicture}
	\begin{scope}[spy using outlines={rectangle,red,magnification=4,size=0.9\textwidth}]
	\node [name=c]{{\includegraphics[height=0.9\textwidth]{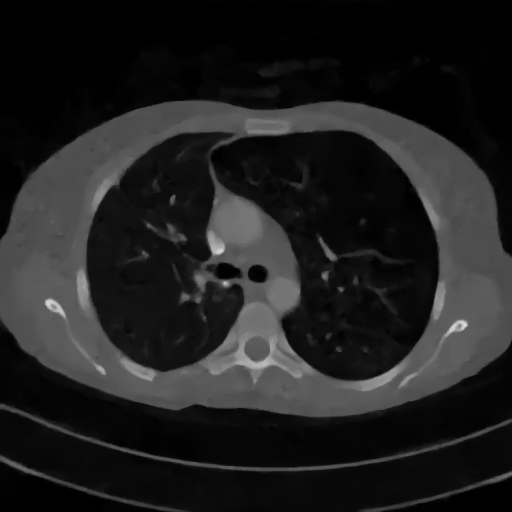}}};
    \spy[magnification=5] on (0.184\textwidth,-0.02\textwidth) in node [name=c1]  at (0,-0.95\textwidth);
	\spy [magnification=4] on (0\textwidth,-0.18\textwidth) in node [name=c1]  at (0,-1.9\textwidth);
	 \spy [magnification=6] on (-0.05\textwidth,0.05\textwidth) in node [name=c1]  at (0,0\textwidth);
    \draw [-stealth, line width=0.1pt, magenta] (-0.02,-0.05) -- ++(-0.1,0.09);
	\end{scope}
	\end{tikzpicture}
	\caption{{ICNN-TV}}
	\label{fig:low_dose-i}
	\end{subfigure}\begin{subfigure}[b]{0.19\textwidth}
	\centering \begin{tikzpicture}
	\begin{scope}[spy using outlines={rectangle,red,magnification=4,size=0.9\textwidth}]
	\node [name=c]{{\includegraphics[height=0.9\textwidth]{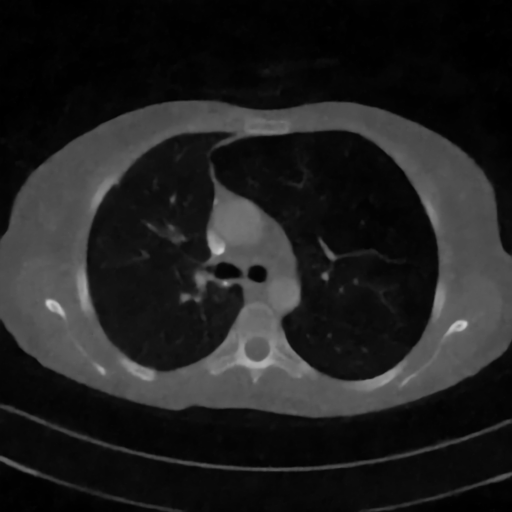}}};
     \spy[magnification=5] on (0.184\textwidth,-0.02\textwidth) in node [name=c1]  at (0,-0.95\textwidth);
	\spy [magnification=4] on (0\textwidth,-0.18\textwidth) in node [name=c1]  at (0,-1.9\textwidth);
	 \spy [magnification=6] on (-0.05\textwidth,0.05\textwidth) in node [name=c1]  at (0,0\textwidth);
    \draw [-stealth, line width=0.1pt, magenta] (-0.02,-0.05) -- ++(-0.1,0.09);
	\end{scope}
	\end{tikzpicture}
	\caption{{GCNN-TV}}
	\label{fig:low_dose-j}
	\end{subfigure}%
	\caption{Three close-ups for each reconstruction by different methods obtained for the chest low-dose CT image. The magenta arrows highlight a region of interest.}
	\label{fig:low_dose}
\end{figure*}

In Figure \ref{fig:low_dose} we report three close-ups of the reconstructions showing different details of the image.
In the first close-up we observe that in some cases the borders of the ascending aorta and superior vena cava sections  pointed by the arrow are not well distinguishable as in the GT image. In particular, we notice that the ICNN and GCNN methods   produce the best images. 
The second crop contains thin  vessels immersed in the dark pulmonary  background. The images obtained with TV, NML and BM3D-WL1 algorithms are too smooth and some  details are hardly visible. In the  BM3D and ICNN-based output images the circular sections of the vessels are distorted into triangular shapes, whereas  the images obtained with gradient-based CNN restore very well the path of the main vessels, without oversmoothing.
In the third row, the close-ups show that ICNN and GCNN recover quite well the object shapes but the former introduces some noisy artifacts, not present in reconstruction obtained from the second which uses CNN gradient-based denoisers. The contribution of the internal TV-based priors in ICNN-TV and GCNN-TV  mitigates the noisy artifacts
(Figures \ref{fig:low_dose-i} and \ref{fig:low_dose-j}).

 We can conclude that in the restoration of this image the use of gradient-based CNN denoiser has some advantages with respect to an image-based CNN denoiser, in better enhancing the objects contours and preserving small details.
We also observe that the BM3D method (Figure \ref{fig:low_dose-e}) always produces  images where the features detachability is enanched, but the anatomic boundaries are unnatural.

Finally, to measure the reconstruction quality and the residual noise, we computed the PSNR and SSIM measures on the whole image and the standard deviation  on a flat the region indicated by  the green square in Figure \ref{fig:low_dose_a}. From the  Table \ref{tab:Deblur_lowdose}, where we report the values obtained for these measures,   we observe that  the GCNN method attains both the best PSNR and SSIM. The BM3D algorithm achieves the second best PSNR but, as we have remarked in the previous comments, it often deformates the curve boundary contours of the objects.
As expected, the TV method over-smooths the image resulting in a very  low standard deviation on the region of interest. In general, adding the  TV as  internal prior in the PnP framework   lowers the standard deviation values of both ICNN and GCNN, as confirmed in the ICNN-TV and GCNN-TV columns.

\begin{table*}
  \begin{center}
 \begin{adjustbox}{max width=\textwidth}
\begin{tabular}{l|c|c|c|c|c|c|c|c}
\thickhline
{} &   TV &    NLM &     BM3D &      BM3D-WL1  & ICNN & GCNN &     ICNN-TV  &    GCNN-TV \\
\thickhline
PSNR & 32.1727 & 30.9899 &  \TextB{34.7675} & 32.9104 &  34.1673 &  \TextBB{35.0309} &  34.0946  &  33.5789  \\
SSIM &   0.9297 &  0.9129 &   \TextB{0.9499} & 0.9358 &   0.9474 &   \TextBB{0.9546} &   0.9466  &   0.9443 \\
ROI-std & \TextGB{0.1746}  &  0.3017 &   0.6569 & 0.5816 &   1.1136 &   1.2366 &   \TextG{0.2844}  &   0.3460 \\
\thickhline
\end{tabular}
\end{adjustbox}
  \end{center}
  \caption{Standard deviation computed on the region of interest inside the green square in Figure \ref{fig:low_dose_a}, for the Low-Dose CT chest images.}
    \label{tab:Deblur_lowdose}
\end{table*}

\section{Conclusions\label{sec:conclusions}}
In this paper we have exploited the flexible Plug-and-Play framework, based on the Half-Quadratic Splitting algorithm, to combine internal and external priors defined on the image gradient domains for the restoration of CT medical images. Moreover we have proposed an external CNN denoiser which acts on the image gradient domain.
 
We have compared the proposed GCNN, ICNN-TV and GCNN-TV  algorithms with  PnP frameworks including other popular denoisers, on synthetic and real CT images characterized by sparse gradients and low-contrast objects.  
The obtained image enhancements confirm that gradient-based priors are effective for the restoration of medical CT images, since the competitors get lower quality indices.
In particular, the proposed GCNN well recovers the curve contours of flat and low-contrast objects, as well as thin vessels, outperforming the ICNN method.
However, the residual artifacts visible on the reconstructions obtained using only an external prior are smoothed out by the ICNN-TV and GCNN-TV algorithms as effect of the internal gradient-based denoiser.   

A deeper insight on the choice of the CNN architecture and the linear feature extractor used into our framework  will be subject of future work, to improve the overall performances of the proposed schemes. 

\section*{Acknowledgments}
This research was funded by the Indam GNCS grant 2020 {\em Ottimizzazione per l'apprendimento automatico e apprendimento automatico per l'ot\-ti\-miz\-za\-zio\-ne}.

\appendix
\renewcommand{\thesection}{\Alph{section}}
\section{Appendix}\label{sec:theorem}
To analyze the convergence properties of Algorithm \ref{alg:hybrid_pnp_hqs}, we start observing that if the denoisers $\mathcal{D}_{\sigma}^{\text{ext}}$ and $\mathcal{D}_{\gamma}^{\text{int}}$ are the proximal maps of two convex functions $g_{1}$ and $g_{2}$, respectively, then   the convergence to a global minimum of the objective function in \eqref{eq:l2+g} is guaranteed \cite{geman1995nonlinear,wang2008new}.
However, in \cite{sreehari2016plug} the authors observe that a denoiser is a proximal map when it is nonexpansive with symmetric gradient, thus limiting the set of suitable denoisers. 
In the effort of allowing less  strict conditions on the involved denoisers, we show in this section that the proposed Algorithm \ref{alg:hybrid_pnp_hqs} satisfies a fixed-point convergence theorem provided only their boundedness.

\begin{definition}[Bounded Denoiser \cite{chan2016plug}] \label{def:bounded_denoiser}
	A \textit{bounded denoiser} is a function $\mathcal{D}_{\epsilon}:\mathbb{R}^l\rightarrow\mathbb{R}^l$ such that for any $\tbold\in\mathbb{R}^l$ and $\epsilon \in \mathbb{R}^{+}$ the following inequality holds:
	\begin{equation}
	\lVert\mathcal{D}_{\epsilon}(\tbold)-\tbold\rVert_2^2\le \epsilon^2 C_{\mathcal{D}}
	\end{equation}
	for a constant $C_{\mathcal{D}}$ independent of $\epsilon$.
\end{definition}

The previous definition entails that given the sequence $(\epsilon_{k})_{k=1}^{+ \infty}$, $\mathcal{D}_{\epsilon_k}$ converges to the identity function of $\mathbb{R}^{l}$ as $\epsilon_k \rightarrow 0$.

In order to state and prove the following  fixed-point theorem, we make some assumptions.\\
Given $(\rho^{\tbold}_k)_{k=1}^{\infty}$ and $(\rho^{\zbold}_k)_{k=1}^{\infty}$ non-decreasing positive sequences, $\Lbold_1 \in \mathbb{R}^{l_1 \times n}$, $\Lbold_2 \in \mathbb{R}^{l_2 \times n}$ as input for Algorithm  \ref{alg:hybrid_pnp_hqs}, then we assume:  %
\\
\begin{enumerate}%
     \item $\mathcal{D}^{\text{ext}}_{\sigma_{k}}$ and $\mathcal{D}^{\text{int}}_{\gamma_{k}}$ are bounded denoisers. \label{ass:A1}
    \item $\Lbold_1$ and $\Lbold_2$ are full-rank matrices. \label{ass:A2}
     \item  $\sum_{k=1}^{+\infty}\sqrt{\dfrac{k}{\rho^{\zbold}_{k}}} < +\infty$, \quad $\sum_{k=1}^{+\infty} \sqrt{\dfrac{k}{\rho^{\tbold}_{k}}} < +\infty$ and $\dfrac{\rho^{\zbold}_{k}}{\rho^{\tbold}_{k}} \to c$ where $c \in \mathbb{R}^{+}$. \label{ass:A3}
\end{enumerate}

\begin{theorem}[Fixed-point convergence theorem for the hybrid PnP algorithm] \label{teo:convergence}
Given the assumptions \ref{ass:A1}-\ref{ass:A3}, there exist $\tbold^* \in \mathbb{R}^{l_1},\zbold^* \in \mathbb{R}^{l_2}$ and $\ubold^* \in \mathbb{R}^{n}$ such that, for   $k\rightarrow\infty$, the following relations hold:

\begin{align}
\tbold_k \to \tbold^*, \quad \Lbold_1\ubold_{k}  \to \tbold^{*}, \quad \zbold_k \to \zbold^*, \quad \Lbold_2\ubold_{k}  \to \zbold^{*}, \quad \ubold_k \to \ubold^*, \nonumber
\end{align}
where $\tbold_k, \zbold_k, \ubold_k $ are computed as in Algorithm \ref{alg:hybrid_pnp_hqs} at step $k$.
\end{theorem}

\begin{proof}

By observing that $\ubold_{k+1}$ is the optimal solution of the minimization problem \eqref{eq:sub_u2}, and by using the relations in \eqref{eq:sigma_gamma} and the assumption \ref{ass:A1},  we get the following chain of inequalities:

\begin{align}
    &\dfrac{1}{2} \lVert \Abold\ubold_{k+1} - \vbold \rVert_{2}^{2} + \dfrac{\rho^{\tbold}_k}{2}\lVert \tbold_{k+1}-\Lbold_1\ubold_{k+1}\rVert_2^2 + \dfrac{\rho^{\zbold}_k}{2}\lVert \zbold_{k+1}-\Lbold_2\ubold_{k+1}\rVert_2^2 \leq \label{eq:LHS_dis}\\
    &\leq \dfrac{1}{2} \lVert \Abold\ubold_{k} - \vbold \rVert_{2}^{2} + \dfrac{\rho^{\tbold}_{k}}{2} \lVert \tbold_{k+1} - \Lbold_1\ubold_{k}\rVert^{2}_{2} + \dfrac{\rho^{\zbold}_k}{2}\lVert \zbold_{k+1}-\Lbold_2\ubold_{k} \lVert_{2}^{2} = \nonumber\\
    & = \dfrac{1}{2} \lVert \Abold\ubold_{k} - \vbold \rVert_{2}^{2} + \dfrac{\rho^{\tbold}_{k}}{2} \lVert \mathcal{D}^{\text{ext}}_{\sigma_{k}}(\Lbold_1\ubold_{k}) - \Lbold_1\ubold_{k}\rVert^{2}_{2} + \dfrac{\rho^{\zbold}_k}{2}\lVert \mathcal{D}^{\text{int}}_{\gamma_{k}}(\Lbold_2\ubold_{k})-\Lbold_2\ubold_{k} \lVert_{2}^{2} \leq \nonumber\\
    & \leq \dfrac{1}{2} \lVert \Abold\ubold_{k} - \vbold \rVert_{2}^{2} + \dfrac{\rho_{k}^{\tbold}}{2} \sigma_{k}^2 C_{\mathcal{D}^{\text{ext}}} + \dfrac{\rho_{k}^{\zbold}}{2} \gamma_{k}^2 C_{\mathcal{D}^{\text{int}}}= \nonumber \\
    & = \dfrac{1}{2} \lVert \Abold\ubold_{k} - \vbold \rVert_{2}^{2} + \dfrac{\alpha}{2} C_{\mathcal{D}^{\text{ext}}} + \dfrac{\beta}{2} C_{\mathcal{D}^{\text{int}}} \leq \nonumber \\
    & = \dfrac{1}{2} \lVert \Abold\ubold_{k} - \vbold \rVert_{2}^{2} + \Tilde{C},\nonumber
\end{align}

with $\Tilde{C}:=\dfrac{\alpha}{2} C_{\mathcal{D}^{\text{ext}}} + \dfrac{\beta}{2} C_{\mathcal{D}^{\text{int}}}$.

Since all the considered terms in \eqref{eq:LHS_dis} are positive, the following inequalities hold: 

\begin{equation}\label{eq:iterative1}
  \dfrac{1}{2} \lVert \Abold\ubold_{k+1} - \vbold \rVert_{2}^{2} \leq  \dfrac{1}{2} \lVert \Abold\ubold_{k} - \vbold \rVert_{2}^{2} + \Tilde{C}
    \leq \dots \leq \dfrac{1}{2} \lVert \Abold\ubold_{1} - \vbold \rVert_{2}^{2} + k\Tilde{C}.
\end{equation}

For the same reason, using \eqref{eq:LHS_dis} and \eqref{eq:iterative1} we get:

\begin{align} 
    & \lVert \tbold_{k+1} - \Lbold_1\ubold_{k+1}\rVert_{2} \leq \sqrt{\dfrac{1}{\rho^{\tbold}_{k}}} \lVert \Abold\ubold_{1} - \vbold \rVert_{2} + \sqrt{\dfrac{2 \Tilde{C} k}{\rho^{\tbold}_k}}, \label{eq:estimate_fond1}\\
    & \lVert \zbold_{k+1} - \Lbold_2\ubold_{k+1}\rVert_{2} \leq \sqrt{\dfrac{1}{\rho^{\zbold}_{k}}} \lVert \Abold\ubold_{1} - \vbold \rVert_{2} + \sqrt{\dfrac{2 \Tilde{C} k}{\rho^{\zbold}_k}}. \label{eq:estimate_fond2}
\end{align}

We now prove that the sequences $(\tbold_{k})_{k=1}^{+ \infty}$ and $(\zbold_{k})_{k=1}^{+ \infty}$ are Cauchy sequences. Starting from the expressions of $\tbold_{k+1}$ and $\zbold_{k+1}$ in Algorithm \ref{alg:hybrid_pnp_hqs}, applying the definition of bounded denoiser and the estimates \eqref{eq:estimate_fond1} and \eqref{eq:estimate_fond2} the following inequalities hold:

\begin{equation}
\begin{split}
    &\lVert \tbold_{k+1}-\tbold_{k}\rVert_2 \le \lVert\mathcal{D}^{\text{ext}}_{\sigma_k}(\Lbold_1\ubold_k)-\Lbold_1\ubold_k\rVert_2 + \lVert \Lbold_1\ubold_k-\tbold_{k}\rVert_2 \le\\ 
    &\leq \sqrt{\dfrac{\alpha}{\rho^{\tbold}_{k}}} \sqrt{C_{\mathcal{D}^{\text{ext}}}} + \sqrt{\dfrac{1}{\rho^{\tbold}_{k-1}}} \lVert \Abold\ubold_{1} - \vbold \rVert_{2} + \sqrt{\dfrac{2 \Tilde{C}(k-1)}{\rho^{\tbold}_{k-1}}}
\end{split}
\end{equation}

\begin{equation}
    \begin{split}
        & \lVert \zbold_{k+1}-\zbold_{k}\rVert_2 \le \lVert \mathcal{D}^{\text{int}}_{\gamma_k}(\Lbold_2\ubold_k)-\Lbold_2\ubold_k\rVert_2 + \lVert \Lbold_2\ubold_k-\zbold_{k}\rVert_2 \le \\
        & \le \sqrt{\dfrac{\beta}{\rho^{\zbold}_{k}}} \sqrt{C_{\mathcal{D}^{\text{int}}}} + \sqrt{\dfrac{1}{\rho^{\zbold}_{k-1}}} \lVert \Abold\ubold_{1} - \vbold \rVert_{2} + \sqrt{\dfrac{2 \Tilde{C}(k-1)}{\rho^{\zbold}_{k-1}}}.
    \end{split}
\end{equation}
By assumption \ref{ass:A3} $(\zbold_{k})_{k=1}^{+ \infty}$ and $(\tbold_{k})_{k=1}^{+ \infty}$ are Cauchy sequences. Hence, there exist $\tbold^{*}$ and $\zbold^{*}$ such that $\tbold_{k} \to \tbold^{*}$ and $\zbold_{k} \to \zbold^{*}$.  

Furthermore,  the following  inequalities (which use  \eqref{eq:estimate_fond1} and \eqref{eq:estimate_fond2}, respectively) state that $\Lbold_1 \ubold_{k+1} \to \tbold^{*}$ and $\Lbold_2 \ubold_{k+1} \to \zbold^{*}$:

\begin{align}
    & \lVert \Lbold_{1} \ubold_{k+1} - \tbold^{*} \rVert_{2} \leq \lVert \Lbold_{1} \ubold_{k+1} - \tbold_{k+1} \rVert_{2} + \lVert \tbold_{k+1} - \tbold^{*} \rVert_{2}, \label{eq:L1u} \\
    & \lVert \Lbold_{2} \ubold_{k+1} - \zbold^{*} \rVert_{2} \leq \lVert \Lbold_{2} \ubold_{k+1} - \zbold_{k+1} \rVert_{2} + \lVert \zbold_{k+1} - \zbold^{*} \rVert_{2} \label{eq:L2u}. 
\end{align}

Now, we prove the convergence of the sequence $(\ubold_{k})_{k=1}^{\infty}$ computed as in Algorithm \ref{alg:hybrid_pnp_hqs}. At step $k$, $\ubold_{k+1}$ is the solution of the convex minimization problem \eqref{eq:sub_u2}, therefore the first order optimality conditions lead: 

\begin{equation}\label{eq:optimality_condition}
    \left(\dfrac{1}{\rho_{k}^{\tbold}}\Abold^{T}\Abold + \Lbold_{1}^T\Lbold_{1} + \dfrac{\rho_{k}^{\zbold}}{\rho_{k}^{\tbold}}\Lbold_{2}^T\Lbold_{2}\right)\ubold_{k+1}=\dfrac{1}{\rho_{k}^{\tbold}}\Abold^T \vbold + \Lbold_{1}^T\tbold_{k+1} + \dfrac{\rho_{k}^{\zbold}}{\rho_{k}^{\tbold}}\Lbold_{2}^T\zbold_{k+1}.
\end{equation}
If we define $\Mbold_{k}:=\frac{1}{\rho_{k}^{\tbold}}\Abold^{T}\Abold + \Lbold_{1}^T\Lbold_{1} + \frac{\rho_{k}^{\zbold}}{\rho_{k}^{\tbold}}\Lbold_{2}^T\Lbold_{2}$, then $\forall\ k > 1$, $\Mbold_{k}$ is invertible for assumption  \ref{ass:A2}.
Hence, we can write for each $k$:
\begin{equation}\label{eq:optimality_condition_1}
   \ubold_{k+1}= \Mbold_{k}^{-1}\left(\dfrac{1}{\rho_{k}^{\tbold}}\Abold^T \vbold + \Lbold_{1}^T\tbold_{k+1} + \dfrac{\rho_{k}^{\zbold}}{\rho_{k}^{\tbold}}\Lbold_{2}^T\zbold_{k+1}\right).
\end{equation}
We observe that the  two sequences in the right hand side of \eqref{eq:optimality_condition_1}, represented by $(\Mbold_{k}^{-1})_{k=1}^{\infty}$ and by the term in parenthesis, are convergent pointwise (by assumption \ref{ass:A3} and by  considering the convergence of the sequences $(\tbold_{k})_{k=1}^{\infty}$ and $(\zbold_{k})_{k=1}^{\infty}$). 
By denoting as $\ubold^{*}$  the product of the two limits, we have proved that $\ubold_{k} \to \ubold^{*}$. \\
This concludes the proof. 
\end{proof}

We point out that this general proof applies also to the algorithm proposed in \cite{zhang2017learning}, for which no convergence results can be found in the literature. Moreover, we believe that with a small effort, our convergence result dealing with multiple denoisers can be extended to ADMM.

The fixed-point convergence Theorem \ref{teo:convergence} entails that the iterations enter in a steady-state and does not guarantee that the fixed-point $\ubold^{*}$ is a minimum of an implicit defined regularized objective as in \eqref{eq:l2+g}. However, in the experimental part, we have shown that the reached fixed-point $\ubold^{*}$ is a very good approximation of the desired image $\ubold$.

\end{document}